\documentclass[reqno,11pt]{article} 
\usepackage[left=2.5cm,right=2.5cm,top=3cm,bottom=3cm,a4paper]{geometry}
\usepackage{amsmath, amssymb}
\usepackage{amsthm, amscd} 
\usepackage{empheq}
\usepackage[all,cmtip]{xy}
\usepackage{float}
\usepackage{caption}
\usepackage{color}
\usepackage{rotating}
\usepackage[utf8]{inputenc}
\usepackage[T1]{fontenc}
\usepackage{url}
\usepackage{booktabs,microtype,needspace}
\usepackage{caption}

\usepackage[titles]{tocloft}
\makeatletter
\DeclareRobustCommand{\rvdots}{%
  \vbox{
    \baselineskip4\p@\lineskiplimit\z@
    \kern-\p@
    \hbox{.}\hbox{.}\hbox{.}
  }}
\makeatother

\newcommand{\Mod}[1]{\ (\textup{mod}\ #1)}
\makeatletter
\def\moverlay{\mathpalette\mov@rlay}
\def\mov@rlay#1#2{\leavevmode\vtop{%
   \baselineskip\z@skip \lineskiplimit-\maxdimen
   \ialign{\hfil$\m@th#1##$\hfil\cr#2\crcr}}}
\newcommand{\charfusion}[3][\mathord]{
    #1{\ifx#1\mathop\vphantom{#2}\fi
        \mathpalette\mov@rlay{#2\cr#3}
      }
    \ifx#1\mathop\expandafter\displaylimits\fi}
\makeatother

\theoremstyle{plain} 
\newtheorem{theorem}{\indent\sc Theorem}[section]
\newtheorem{lemma}[theorem]{\indent\sc Lemma}

\newtheorem{proposition}[theorem]{\indent\sc Proposition}

\theoremstyle{definition} 
\newtheorem{definition}[theorem]{\indent\sc Definition}

\newtheorem{remark}[theorem]{\indent\sc Remark}

\newtheorem{thmx}{Theorem}

\makeatletter
\def\address#1#2{\begingroup
\noindent\parbox[t]{7.8cm}{%
\small{\scshape\ignorespaces#1}\par\vskip1ex
\noindent\small{\itshape E-mail address}%
\/: #2\par\vskip4ex}\hfill%
\endgroup}%
\makeatother

\title{On the Hasse-Sugawara problem and its application to Hilbert's twelfth problem}
\author{
\textsc{Ho Yun Jung$^*$, Ja Kyung Koo and Dong Hwa Shin} 
}
\date{} 
\begin{document}

\allowdisplaybreaks

\maketitle

\footnote{
$^*$Corresponding author.\\
\thanks{
}
}

\begin{abstract}
We give an analytic proof of the Hasse-Sugawara problem 
by treating the seven
remaining families of moduli. 
We also construct an explicit integer translate of a Weber value
whose relative norms generate every intermediate field, and
express these generators using extended form class groups. 
\end{abstract}

\tableofcontents

\section {Introduction}

The generation of abelian extensions of imaginary quadratic fields
by special values of modular functions is a classical problem in complex multiplication. 
We study the generation of ray class fields by a single Weber value
and of their intermediate fields by relative norms.
\par
Let $K$ be an imaginary quadratic field of discriminant $d_K$, and let $\mathcal{O}_K$ be its ring of integers.
For a nonzero ideal $\mathfrak{m}$ of $\mathcal{O}_K$,
let $C(\mathfrak{m})$ denote the ray class group of $K$ modulo $\mathfrak{m}$.  
By class field theory,
there is a unique finite abelian extension $K_\mathfrak{m}$ of $K$,
unramified outside the prime ideals of $\mathcal{O}_K$ dividing $\mathfrak{m}$,
for which the Artin map modulo $\mathfrak{m}$ induces an isomorphism
\begin{equation*}
C(\mathfrak{m})\stackrel{\sim}
{\rightarrow}\mathrm{Gal}(K_\mathfrak{m}/K).
\end{equation*}
This extension is called the 
ray class field of $K$ modulo $\mathfrak{m}$; 
see \cite[Chapters III--IV]{Janusz}.
When $\mathfrak{m}=\mathcal{O}_K$, it is the Hilbert class field $H_K$ of $K$. 
\par
Let $E/\mathbb{C}$ be an elliptic curve associated with the lattice $\mathcal{O}_K$. Then
$E$ has complex multiplication by $\mathcal{O}_K$, and scalar multiplication gives an isomorphism
\begin{equation*}
\mathcal{O}_K\stackrel{\sim}{\rightarrow}
\mathrm{End}(E),\quad a\mapsto[a]_E.
\end{equation*}
In particular, the group $\mathrm{Aut}(E)$ of automorphisms 
is isomorphic to $\mathcal{O}_K^\times$. 
We use the Weierstrass model
\begin{equation}\label{model}
y^2=4x^3-g_2x-g_3,\quad\textrm{where}~g_2=g_2(\mathcal{O}_K)~
\textrm{and}~g_3=g_3(\mathcal{O}_K)
\end{equation}
for $E$. 
The associated complex analytic isomorphism is
\begin{equation}\label{parametrization}
\phi_E:\mathbb{C}/\mathcal{O}_K\stackrel{\sim}{\rightarrow}
E(\mathbb{C})~(\subset\mathbb{P}^2(\mathbb{C})),
\quad z+\mathcal{O}_K\mapsto\left\{
\begin{array}{ll}
\left[\wp(z;\,\mathcal{O}_K):\wp'(z;\,\mathcal{O}_K):1\right] & \textrm{if}~z\in\mathbb{C}\setminus\mathcal{O}_K,\\
\mathrm{[}0:1:0] & \textrm{if}~z\in\mathcal{O}_K.
\end{array}\right.
\end{equation}
Here $\wp=\wp(\,\cdot\,;\,\mathcal{O}_K):\mathbb{C}\rightarrow
\mathbb{C}\cup\{\infty\}$
is the Weierstrass function for the lattice $\mathcal{O}_K$,
and $\wp'$ is its derivative;
see
\cite[Theorem 10.14]{Cox} and 
\cite[Proposition 3.6 in Chapter VI]{Silverman09}.
Let $O_E=[0:1:0]$ denote the origin of $E$. 
\par
We consider the Weber function 
\begin{equation*}
h_E:E\rightarrow\mathbb{P}^1(\mathbb{C})=\mathbb{C}\cup\{\infty\}
\end{equation*}
defined by $h_E(O_E)=\infty$ and
\begin{equation}\label{Weberh}
h_E(\phi_E(z+\mathcal{O}_K))=\left\{\begin{array}{ll}
\displaystyle\frac{g_2g_3}{\Delta}\,\wp(z;\,\mathcal{O}_K) &
\textrm{if}~ j(E)\neq1728,\,0,\vspace{0.1cm}\\
\displaystyle\frac{g_2^2}{\Delta}\,\wp(z;\,\mathcal{O}_K)^2 & \textrm{if}~j(E)=1728,\vspace{0.1cm}\\
\displaystyle\frac{g_3}{\Delta}\,\wp(z;\,\mathcal{O}_K)^3 & \textrm{if}~j(E)=0
\end{array}\right.\quad(z\in\mathbb{C}\setminus\mathcal{O}_K),
\end{equation}
where $\Delta=g_2^3-27g_3^2$. 
It realizes the quotient map 
$E\rightarrow E/\mathrm{Aut}(E)\simeq\mathbb{P}^1(\mathbb{C})$.
See \cite[Example 5.5.2]{Silverman94}.
\par
For a nonzero proper ideal $\mathfrak{m}$ of $\mathcal{O}_K$,
a point $P\in E(\mathbb{C})$ is called a \textit{primitive 
$\mathfrak{m}$-division point} if 
\begin{equation*}
\mathrm{Ann}_{\mathcal{O}_K}(P)=\{a\in\mathcal{O}_K\mid
[a]_EP=O_E\}=\mathfrak{m}.
\end{equation*}
The classical theory of complex multiplication, developed
by Kronecker, Weber, Hasse, Deuring, and others, gives
\begin{equation*}
H_K=K(j(E))
\end{equation*}
and
\begin{equation}\label{overH}
K_\mathfrak{m}=H_K(h_E(P))\quad\textrm{for any primitive $\mathfrak{m}$-division point}~P\in 
E(\mathbb{C}).
\end{equation}
As $P$ varies over these points, the values 
$h_E(P)$ form the full set of 
Galois conjugates over $H_K$ of any one of them; 
see \cite[$\S$10.3]{Lang}. 
Hasse further asked whether $H_K$ in \eqref{overH} can be replaced by 
$K$, that is, whether
\begin{equation}\label{overK}
K_\mathfrak{m}=K(h_E(P))\quad\textrm{for every primitive
$\mathfrak{m}$-division point}~P\in E(\mathbb{C}).
\end{equation}
\par
For a nonzero ideal $\mathfrak{a}$ of $\mathcal{O}_K$, let 
$\mathrm{N}(\mathfrak{a})=|\mathcal{O}_K/\mathfrak{a}|$ be its absolute norm. 
Unless otherwise specified, we denote by $\mathfrak{p}$ a nonzero prime ideal 
of $\mathcal{O}_K$. 
Put
\begin{equation*}
\varphi(\mathfrak{a})=\mathrm{N}(\mathfrak{a})\prod_{\mathfrak{p}\mid\mathfrak{a}}
\left(1-\frac{1}{\mathrm{N}(\mathfrak{p})}\right)\quad
\textrm{and}\quad\Psi(\mathfrak{a})=\mathrm{N}(\mathfrak{a})\prod_{\mathfrak{p}|\mathfrak{a}}
\left(1-\frac{2}{\mathrm{N}(\mathfrak{p})}\right). 
\end{equation*}
In 1936 Sugawara \cite{Sugawara1, Sugawara2} gave
the first partial answers to Hasse's question, proving \eqref{overK}
whenever at least one of the following 
conditions holds:
\begin{enumerate}
\item[(S1)] $4\mathcal{O}_K\mid\mathfrak{m}$ and $\mathfrak{m}\neq4\mathcal{O}_K$;
\item[(S2)] $\mathfrak{p}^2\mid\mathfrak{m}$
for some prime ideal $\mathfrak{p}$ of $\mathcal{O}_K$ such that $\mathfrak{p}\mid2\mathcal{O}_K$, 
and $\varphi(\mathfrak{m})\geq6$;
\item[(S3)] $\Psi(\mathfrak{m})\geq5$. 
\end{enumerate}
His method uses addition formulas for the Weierstrass $\wp$-function to recover $j(E)$ from suitable Galois conjugates of $h_E(P)$.
Subsequently, 
Jung-Koo-Shin \cite{J-K-S} and Koo-Shin-Yoon 
\cite{K-S-Y} proved \eqref{overK} for $\mathfrak{m}=N\mathcal{O}_K$,
where $N\geq2$ is an integer, under the hypothesis that
$\mathfrak{m}$ is the exact conductor of $K_\mathfrak{m}/K$,
using modular functions and the second Kronecker limit formula.  
\par
Let $\mathfrak{f}_{K_\mathfrak{m}/K}$ denote the exact conductor of $K_\mathfrak{m}/K$. Since
\begin{equation*}
\mathfrak{f}_{K_\mathfrak{m}/K}\mid\mathfrak{m}\quad\textrm{and}
\quad
K_\mathfrak{m}=K_{\mathfrak{f}_{K_\mathfrak{m}/K}}, 
\end{equation*}
it is natural to consider \eqref{overK}
under the hypothesis $\mathfrak{m}=\mathfrak{f}_{K_\mathfrak{m}/K}$.
Morton \cite{Morton} recently
presented a proof of \eqref{overK} under this hypothesis. 
Starting from Sugawara's criteria
and the conductor condition, he reduced the problem
to seven remaining families of moduli and
treated them by algebraic and arithmetic arguments
involving explicit calculations on 
suitable elliptic curves. 
\par
In this paper, we give an alternative proof of \eqref{overK}. 
We follow Morton's reduction to the seven remaining families,
recalled in Proposition \ref{residual}, 
and treat them by two methods. 
Four families are handled by a character argument involving 
Siegel functions, the second Kronecker limit formula and 
the nonvanishing of Hecke $L$-values at $s=1$ (Theorem \ref{main1}). 
The other three are handled by explicit
estimates for the Weierstrass $\wp$-function, the elliptic modular function $j$
and the modular discriminant $\Delta$ (Theorem \ref{main2}). 
Our first main result is the following. 

\begin{thmx}[Theorem \ref{conclusion}]\label{main}
Let $K$ be an imaginary quadratic field,
and let $E$ and $h_E$ be as in \eqref{model}--\eqref{Weberh}. 
If $\mathfrak{m}$ is a nonzero proper ideal of $\mathcal{O}_K$
that is the exact conductor of $K_\mathfrak{m}/K$, then
\begin{equation*}
K_\mathfrak{m}=K(h_E(P))
\end{equation*}
for every primitive $\mathfrak{m}$-division
point $P\in E(\mathbb{C})$.
\end{thmx}

Ramachandra \cite[Theorems 9 and 10]{Ramachandra} constructed ray class units whose relative norms, 
as well as all positive integral powers of these norms, generate every 
finite abelian intermediate field over $K$. 
He also observed that this norm-generation property can fail for simpler invariants \cite[pp. 132--133]{Ramachandra}. 
Schertz \cite[Conjecture 6.8.3]{Schertz} formulated related conjectures on the generation of intermediate fields by relative
norms of simpler elliptic invariants. In this paper,
however,  
we obtain the same norm-generation property by subtracting an explicit integer
from a Weber value (Theorem \ref{intermediate}). This gives the following explicit formulation of
the solution to Hilbert's twelfth problem for imaginary quadratic fields. 

\begin{thmx}[Theorem \ref{Hilbert}]
Let $K$ be an imaginary quadratic field of discriminant $d_K$, 
and let $L$ be a finite extension of $K$. 
Then $L/K$ is abelian if and only if 
there exist a positive integer $N$ divisible by $4$ and a subgroup $H$ of 
the extended form class group $C_{\Gamma_1(N)}(d_K,\,N)$
such that
\begin{equation*}
L=K\left(
\prod_{[Q]_{\Gamma_1(N)}\in H}\left(f_{\left[\begin{smallmatrix}
0 &a_Q'/N
\end{smallmatrix}\right]}^{(|\mathcal{O}_K^\times|/2)}(-\overline{\omega}_Q)
-\mathcal{C}_{N\mathcal{O}_K}\right)\right).
\end{equation*}
\end{thmx}

Here, $\omega_Q\in\mathbb{H}$ is the zero of $Q(x,\,1)$,
and $a_Q'$ is an inverse modulo $N$ of the leading coefficient of $Q$. 
The functions $f_\mathbf{v}^{(k)}$, the integers $\mathcal{C}_\mathfrak{m}$ and
the extended form class groups $C_{\Gamma_1(N)}(d_K,\,N)$
are defined in Sections \ref{sect3}, \ref{sect8} and \ref{sect9}. 

\bigskip

Throughout this paper, we use the following notation unless otherwise specified: 

\begin{itemize}
\item $K$ : an imaginary quadratic field of discriminant $d_K$.
\item $\mathcal{O}_K$ : the ring of integers of $K$.
\item $\mathfrak{m}$: a nonzero proper ideal of $\mathcal{O}_K$.
\item $\mathfrak{p}$ : a nonzero prime ideal of $\mathcal{O}_K$.  
\item $E$ : the elliptic curve defined by the Weierstrass model in \eqref{model}
with the parametrization $\phi_E:\mathbb{C}/\mathcal{O}_K
\stackrel{\sim}{\rightarrow}E(\mathbb{C})$
in \eqref{parametrization}.
\end{itemize}
Furthermore, we write Galois actions on the right, denoting the action of 
$\sigma$ on $a$ by $a^\sigma$. 

\section {Remaining families of moduli}

Following Morton \cite[Section 1]{Morton}, we reduce
the proof of Theorem \ref{main} to seven families of moduli
not covered by the earlier results.  For the sake of completeness, 
we give a proof of the reduction. 
\par
Let $h_K$ denote the class number of $K$. By the results of Heegner \cite{Heegner}, Baker \cite{Baker1, Baker2, Baker3} 
and Stark
\cite{Stark}, we have
\begin{equation}\label{classnumberone}
h_K=1\quad\Longleftrightarrow\quad
d_K\in\{-3,\,-4,\,-7,\,-8,\,-11,\,-19,\,-43,\,-67,\,-163\} 
\end{equation}
(see also \cite[Theorem 12.34]{Cox}).
For a positive integer $N$, let $\zeta_N=e^{2\pi\mathrm{i}/N}$. 
We note that
\begin{equation}\label{O_K}
\mathcal{O}_K=[\tau_K,\,1]
\quad\textrm{where}~
\tau_K=\left\{\begin{array}{ll}
\displaystyle\frac{\sqrt{d_K}}{2} & \textrm{if}~d_K\equiv0\Mod{4},\vspace{0.1cm}\\
\displaystyle\frac{-1+\sqrt{d_K}}{2} & \textrm{if}~d_K\equiv1\Mod{4},
\end{array}\right.
\end{equation}
and 
\begin{equation}\label{O_Ktimes}
\mathcal{O}_K^\times=\left\{\begin{array}{ll}
\langle\zeta_4\rangle& \textrm{if}~K=\mathbb{Q}(\sqrt{-1})~(\textrm{equivalently,}~d_K=-4),\\
\langle\zeta_6\rangle & \textrm{if}~K=\mathbb{Q}(\sqrt{-3})~(\textrm{equivalently,}~d_K=-3),\\
\langle\zeta_2\rangle & \textrm{otherwise} 
\end{array}\right.
\end{equation}
(see \cite[(5.13) and p. 92]{Cox}). 
\par
For a nonzero ideal $\mathfrak{m}$ of $\mathcal{O}_K$, the ray class group
$C(\mathfrak{m})$ is defined by 
\begin{equation*}
C(\mathfrak{m})=I_K(\mathfrak{m})/P_{K,\,1}(\mathfrak{m}),
\end{equation*}
where $I_K(\mathfrak{m})$ is the group of fractional ideals of $K$ relatively prime to $\mathfrak{m}$ and
$P_{K,\,1}(\mathfrak{m})$ is its subgroup given by
\begin{equation*}
P_{K,\,1}(\mathfrak{m})=\langle\nu\mathcal{O}_K\mid\nu\in\mathcal{O}_K\setminus\{0\},~\nu\equiv1\Mod{\mathfrak{m}}\rangle. 
\end{equation*}
Let $K_\mathfrak{m}$ denote the corresponding ray class field so that
$C(\mathfrak{m})\simeq\mathrm{Gal}(K_\mathfrak{m}/K)$ via
the Artin map modulo $\mathfrak{m}$. 
Note that
if $\mathfrak{a}$ and $\mathfrak{b}$ are nonzero ideals of $\mathcal{O}_K$ such that
$\mathfrak{a}\mid\mathfrak{b}$, then $K_\mathfrak{a}\subseteq K_\mathfrak{b}$
(see \cite[Theorem 9.9 in Chapter V]{Janusz}). 
The conductor of the extension $K_\mathfrak{m}/K$
is defined as the smallest (with respect to divisibility) nonzero ideal $\mathfrak{f}=\mathfrak{f}_{K_\mathfrak{m}/K}$ of $\mathcal{O}_K$ such that $K_\mathfrak{m}\subseteq K_\mathfrak{f}$,
that is, there is no nonzero ideal $\mathfrak{a}$ of $\mathcal{O}_K$ such that
\begin{equation*}
\mathfrak{a}\mid\mathfrak{f},\quad\mathfrak{a}\neq\mathfrak{f},\quad K_\mathfrak{m}\subseteq K_\mathfrak{a}. 
\end{equation*}
\par
For a rational prime $p$, we write the prime ideal factorization of $p\mathcal{O}_K$
as follows:
\begin{align*}
p\mathcal{O}_K=\left\{
\begin{array}{ll}
\mathfrak{p}_{p,\,\mathrm{split}}
\mathfrak{p}'_{p,\,\mathrm{split}}
 & \textrm{if $p$ splits in $K$},\\
\mathfrak{p}_{p,\,\mathrm{inert}}  & \textrm{if $p$ is inert in $K$},\\
\mathfrak{p}_{p,\,\mathrm{ramified}}^2  & \textrm{if $p$ ramifies in $K$}.
\end{array}\right.
\end{align*}
When $p$ is not inert in $K$, we write
$\mathfrak{p}_{p,\,\mathrm{noninert}}$ 
for a prime ideal of $\mathcal{O}_K$ lying above $p$. 

\begin{proposition}\label{residual}
Let $K$ be an imaginary quadratic field such that $h_K>1$,
and let $\mathfrak{m}$ be a nonzero proper ideal of $\mathcal{O}_K$ such that
$\mathfrak{f}_{K_\mathfrak{m}/K}=\mathfrak{m}$. 
Assume that $\mathfrak{m}$ satisfies none of \textup{(S1), (S2), (S3)} and that
\begin{equation}\label{Nassumption}
\mathfrak{m}\neq N\mathcal{O}_K\quad\textrm{for every integer $N>1$}.
\end{equation} 
Then, up to interchanging conjugate prime ideals,
$\mathfrak{m}$ has one of the following seven forms:
\begin{equation}
\begin{gathered}
\mathfrak{p}_{2,\,\mathrm{split}}^2\mathfrak{p}_{3,\,\mathrm{noninert}},\quad
\mathfrak{p}_{2,\,\mathrm{split}}^3,\quad
\mathfrak{p}_{3,\,\mathrm{split}}^2,\quad
\mathfrak{p}_{3,\,\mathrm{split}}^2\mathfrak{p}'_{3,\,\mathrm{split}},\\
\mathfrak{p}_{5,\,\mathrm{noninert}},\quad
\mathfrak{p}_{3,\,\mathrm{noninert}}\mathfrak{p}_{5,\,\mathrm{noninert}},\quad
\mathfrak{p}_{3,\,\mathrm{split}}\mathfrak{p}'_{3,\,\mathrm{split}}\mathfrak{p}_{5,\,\mathrm{noninert}}.
\end{gathered}
\end{equation}
\end{proposition}
\begin{proof}
Since $h_K>1$, we have $K\neq\mathbb{Q}(\sqrt{-1}),\,\mathbb{Q}(\sqrt{-3})$ by \eqref{classnumberone},
and hence $\mathcal{O}_K^\times=\{\pm1\}$ by \eqref{O_Ktimes}. 
For a nonzero ideal $\mathfrak{a}$ of $\mathcal{O}_K$, let
\begin{equation*}
U_1(\mathfrak{a}):=\{u\in\mathcal{O}_K^\times\mid u\equiv1\Mod{\mathfrak{a}}\}, 
\end{equation*}
and then
\begin{equation}\label{|U|}
\quad|U_1(\mathfrak{a})|=\left\{\begin{array}{ll}
2 & \textrm{if}~\mathfrak{a}\mid2\mathcal{O}_K,\\
1 & \textrm{if}~\mathfrak{a}\nmid2\mathcal{O}_K.
\end{array}\right.
\end{equation}
The exact sequence
\begin{equation}\label{exact}
1\rightarrow U_1(\mathfrak{a})\rightarrow
\mathcal{O}_K^\times\rightarrow
(\mathcal{O}_K/\mathfrak{a})^\times
\rightarrow C(\mathfrak{a})\rightarrow C(\mathcal{O}_K)
\rightarrow 1
\end{equation}
gives the extension degree formula
\begin{equation*}
[K_\mathfrak{a}:H_K]=\frac{\varphi(\mathfrak{a})|U_1(\mathfrak{a})|}{2}
\end{equation*}
(see \cite[Theorem 1.7 in Chapter V]{Milne}). It follows that
if $\mathfrak{p}\mid\mathfrak{a}$ and $\mathfrak{p}^e\|\mathfrak{a}$, then 
\begin{equation}\label{degreeformula}
[K_\mathfrak{a}:K_{\mathfrak{a}\mathfrak{p}^{-1}}]=
\frac{\lambda_\mathfrak{p}(e)}{[U_1(\mathfrak{a}\mathfrak{p}^{-1}):
U_1(\mathfrak{a})]}
\quad\textrm{where}~\lambda_\mathfrak{p}(e):=\left\{\begin{array}{ll}
\mathrm{N}(\mathfrak{p})-1 & \textrm{if}~e=1,\\
\mathrm{N}(\mathfrak{p}) & \textrm{if}~e\geq2.  
\end{array}\right.
\end{equation}
Thus the hypothesis $\mathfrak{f}_{K_\mathfrak{m}/K}=\mathfrak{m}$ is equivalent to 
\begin{align}\label{exactcondition}
\frac{\lambda_\mathfrak{p}(e_\mathfrak{p})}{[U_1(\mathfrak{m}\mathfrak{p}^{-1}):
U_1(\mathfrak{m})]}=[K_\mathfrak{m}:K_{\mathfrak{m}\mathfrak{p}^{-1}}]>1~
\textrm{for every}~\mathfrak{p}\mid\mathfrak{m}~\textrm{such that}~\mathfrak{p}^{e_\mathfrak{p}}
\|\mathfrak{m}.
\end{align}
There are three possible cases:
\begin{equation*}
\left\{\begin{array}{l}
\gcd(\mathfrak{m},\,2\mathcal{O}_K)\neq\mathcal{O}_K~\textrm{and 
$2$ is not inert in $K$,\quad or}\\
\gcd(\mathfrak{m},\,2\mathcal{O}_K)\neq\mathcal{O}_K~\textrm{and 
$2$ is inert in $K$,\quad or}\\
\gcd(\mathfrak{m},\,2\mathcal{O}_K)=\mathcal{O}_K.
\end{array}\right.
\end{equation*}
\begin{enumerate}
\item[Case 1.] Suppose that
$\gcd(\mathfrak{m},\,2\mathcal{O}_K)\neq\mathcal{O}_K$ and
$2$ is not inert in $K$. Let $\mathfrak{p}\mid
\gcd(\mathfrak{m},\,2\mathcal{O}_K)$, and write
$\mathfrak{m}=\mathfrak{p}^e\mathfrak{b}$ for an ideal $\mathfrak{b}$ of 
$\mathcal{O}_K$ such that
$\gcd(\mathfrak{b},\,\mathfrak{p})=\mathcal{O}_K$. 
Since $\mathrm{N}(\mathfrak{p})=2$ and $\mathfrak{f}_{K_\mathfrak{m}/K}=\mathfrak{m}$, 
we deduce by \eqref{exactcondition} that $e\geq2$. 
Since (S2) fails, we must have
\begin{equation}\label{6}
6>\varphi(\mathfrak{m})=2^{e-1}\varphi(\mathfrak{b}). 
\end{equation}
Hence we get $e\in\{2,\,3\}$.
\begin{enumerate}
\item[Case 1-1.] Let $e=2$, so $\varphi(\mathfrak{b})<3$. 
The possibilities
for $\mathfrak{b}$ are 
\begin{equation*}
\mathfrak{b}=\left\{\begin{array}{ll}
\mathcal{O}_K,\quad\mathfrak{p}'_{2,\,\mathrm{split}},
\quad\mathfrak{p}^{'2}_{2,\,\mathrm{split}},
\quad\mathfrak{p}'_{2,\,\mathrm{split}}
\mathfrak{p}_{3,\,\mathrm{noninert}},
\quad\mathfrak{p}_{3,\,\mathrm{noninert}}
 & \textrm{if}~\mathfrak{p}=\mathfrak{p}_{2,\,\mathrm{split}},\\
\mathcal{O}_K,\quad
\mathfrak{p}_{3,\,\mathrm{noninert}} & \textrm{if}~\mathfrak{p}=\mathfrak{p}_{2,\,\mathrm{ramified}}.
\end{array}\right.
\end{equation*}
Furthermore, using \eqref{Nassumption} and \eqref{exactcondition}, we find that
\begin{equation}\label{Case2-1}
\mathfrak{m}=\mathfrak{p}_{2,\,\mathrm{split}}^2\mathfrak{p}_{3,\,\mathrm{noninert}}. 
\end{equation}
\item[Case 1-2.] Let $e=3$. Then \eqref{6} yields $\varphi(\mathfrak{b})<\frac{3}{2}$, 
and hence
\begin{equation*}
\mathfrak{b}=\left\{\begin{array}{ll}
\mathcal{O}_K,\quad\mathfrak{p}'_{2,\,\mathrm{split}}
 & \textrm{if}~\mathfrak{p}=\mathfrak{p}_{2,\,\mathrm{split}},\\
\mathcal{O}_K & \textrm{if}~\mathfrak{p}=\mathfrak{p}_{2,\,\mathrm{ramified}}.
\end{array}\right.
\end{equation*}
Applying \eqref{exactcondition}, we obtain
\begin{equation}\label{Case2-2}
\mathfrak{m}=\mathfrak{p}_{2,\,\mathrm{split}}^3.
\end{equation}
\end{enumerate}
\item[Case 2.] Suppose that $\gcd(\mathfrak{m},\,2\mathcal{O}_K)\neq
\mathcal{O}_K$ and $2$ is inert in $K$. Let $\mathfrak{p}=\mathfrak{p}_{2,\,\mathrm{inert}}$
and so $\mathrm{N}(\mathfrak{p})=4$.  
Write $\mathfrak{m}=\mathfrak{p}^e\mathfrak{b}$ for an ideal $\mathfrak{b}$ of $\mathcal{O}_K$ such that $\gcd(\mathfrak{b},\,\mathfrak{p})=
\mathcal{O}_K$. If $e\geq2$, then $4\mathcal{O}_K\mid\mathfrak{m}$ and then
$\mathfrak{m}=4\mathcal{O}_K$ since (S1) fails. But this contradicts \eqref{Nassumption}, and
hence $e=1$. Since (S3) fails, we have
\begin{equation}\label{5a}
5>\Psi(\mathfrak{m})=\Psi(\mathfrak{p}\mathfrak{b})=2\Psi(\mathfrak{b}). 
\end{equation}
For a prime ideal $\mathfrak{q}$ of $\mathcal{O}_K$ of odd norm
and an integer $r\geq1$, we see that
\begin{equation*}
\Psi(\mathfrak{q}^r)=\mathrm{N}(\mathfrak{q})^{r-1}(\mathrm{N}(\mathfrak{q})-2)
\quad\textrm{and}\quad \mathrm{N}(\mathfrak{q})\geq3. 
\end{equation*} 
And it follows from \eqref{5a} that
\begin{equation*}
\mathfrak{b}=\mathcal{O}_K,\quad\mathfrak{p}_{3,\,\mathrm{noninert}},\quad
\mathfrak{p}_{3,\,\mathrm{split}}\mathfrak{p}'_{3,\,\mathrm{split}}, 
\end{equation*}
and the assumption \eqref{Nassumption} forces us to get
\begin{equation*}
\mathfrak{m}=\mathfrak{p}_{2,\,\mathrm{inert}}\mathfrak{p}_{3,\,\mathrm{noninert}}. 
\end{equation*}
However, we observe by \eqref{|U|} that
\begin{equation*}
\frac{\lambda_{\mathfrak{p}_{3,\,\mathrm{noninert}}}(1)}{[U_1(\mathfrak{p}_{2,\,\mathrm{inert}}):U_1(\mathfrak{m})]}=\frac{3-1}{2}=1,
\end{equation*}
contradicting \eqref{exactcondition}. 
Thus Case 2 yields no remaining form. 
\item[Case 3.] Suppose that $\gcd(\mathfrak{m},\,2\mathcal{O}_K)=\mathcal{O}_K$.
Let $\mathfrak{m}=\prod_{\mathfrak{p}\mid\mathfrak{m}}\mathfrak{p}^{e_\mathfrak{p}}$. 
Since (S3) fails, we attain
\begin{equation*}
5>\Psi(\mathfrak{m})=\prod_{\mathfrak{p}\mid\mathfrak{m}}
\Psi(\mathfrak{p}^{e_\mathfrak{p}})=
\prod_{\mathfrak{p}\mid\mathfrak{m}}\mathrm{N}(\mathfrak{p})^{e_\mathfrak{p}-1}(\mathrm{N}(\mathfrak{p})-2)\quad\textrm{and}
\quad
\mathrm{N}(\mathfrak{p})\geq3~(\mathfrak{p}\mid\mathfrak{m}). 
\end{equation*} 
Therefore we find that
\begin{align*}
\mathfrak{m}=\left\{\begin{array}{l}\mathfrak{p}_{3,\,\mathrm{noninert}},\quad
\mathfrak{p}_{3,\,\mathrm{noninert}}^2,\quad
\mathfrak{p}_{3,\,\mathrm{split}}\mathfrak{p}'_{3,\,\mathrm{split}},\quad
\mathfrak{p}_{3,\,\mathrm{split}}^2\mathfrak{p}'_{3,\,\mathrm{split}},\\
\mathfrak{p}_{5,\,\mathrm{noninert}},\quad
\mathfrak{p}_{3,\,\mathrm{noninert}}\mathfrak{p}_{5,\,\mathrm{noninert}},\quad
\mathfrak{p}_{3,\,\mathrm{split}}\mathfrak{p}'_{3,\,\mathrm{split}}\mathfrak{p}_{5,\,\mathrm{noninert}}.
\end{array}\right.
\end{align*}
Applying \eqref{Nassumption} and \eqref{exactcondition}, we derive
the following five forms:
\begin{equation}\label{Case3}
\mathfrak{m}=
\mathfrak{p}_{3,\,\mathrm{split}}^2,\quad
\mathfrak{p}_{3,\,\mathrm{split}}^2\mathfrak{p}'_{3,\,\mathrm{split}},\quad
\mathfrak{p}_{5,\,\mathrm{noninert}},\quad
\mathfrak{p}_{3,\,\mathrm{noninert}}\mathfrak{p}_{5,\,\mathrm{noninert}},\quad
\mathfrak{p}_{3,\,\mathrm{split}}\mathfrak{p}'_{3,\,\mathrm{split}}\mathfrak{p}_{5,\,\mathrm{noninert}}.
\end{equation}
Combining \eqref{Case2-1}, \eqref{Case2-2} and \eqref{Case3} completes the proof. 
\end{enumerate}
\end{proof}

\section {Fields of modular functions}\label{sect3}

We recall the modular functions in order to express Weber values
and the Galois action on them. Their Fourier expansions and 
the identity relating differences of Fricke functions to Siegel functions
provide the analytic tools for the proofs below. 
\par
Let $\mathbb{H}=\{\tau\in\mathbb{C}\mid\mathrm{Im}(\tau)>0\}$ be the complex upper 
half-plane. For $z_1,\,z_2\in\mathbb{C}$ such that $z_1/z_2\in\mathbb{H}$, 
we denote by $[z_1,\,z_2]$ the lattice $\mathbb{Z}z_1+\mathbb{Z}z_2$.  
Define the functions $g_2$, $g_3$ and $\Delta$ on $\mathbb{H}$ by 
\begin{equation}\label{Deltadef}
g_2(\tau)=g_2([\tau,\,1]),\quad
g_3(\tau)=g_3([\tau,\,1]),\quad
\Delta(\tau)=\Delta([\tau,\,1])=g_2(\tau)^3-27g_3(\tau)^2\quad(\tau\in\mathbb{H}), 
\end{equation}
which are modular forms for $\mathrm{SL}_2(\mathbb{Z})$
of weight $4$, $6$, $12$, respectively (see \cite[pp. 36--38]{Lang}).
Put
\begin{equation*}
q_\tau=e^{2\pi\mathrm{i}\tau}\quad(\tau\in\mathbb{H}).
\end{equation*}

\begin{lemma}\label{g2g3Deltaexpansion}
The functions $g_2$, $g_3$, $\Delta$ satisfy the following expansion formulas:
\begin{align}
\frac{12}{(2\pi\mathrm{i})^4}\,g_2(\tau)&=1+240\sum_{n=1}^\infty\frac{n^3q_\tau^n}{1-q_\tau^n}\qquad(\tau\in\mathbb{H}),\label{g2expansion}\\
\frac{6^3}{(2\pi\mathrm{i})^6}\,g_3(\tau)&=-1+504\sum_{n=1}^\infty\frac{n^5q_\tau^n}{1-q_\tau^n}
\qquad(\tau\in\mathbb{H}),\label{g3expansion}\\
\frac{1}{(2\pi\mathrm{i})^{12}}\,\Delta(\tau)&=q_\tau\prod_{n=1}^\infty (1-q_\tau^n)^{24}
\label{Deltaexpansion}
\qquad(\tau\in\mathbb{H}).
\end{align}
\end{lemma}
\begin{proof}
See \cite[Proposition 4 in Chapter 4 and Theorem 5 in Chapter 18]{Lang}.
\end{proof}

For later use, we also introduce the well-known expansion formula of the Weierstrass function. 
For $z\in\mathbb{C}$, we also use the notation 
$q_z=e^{2\pi\mathrm{i}z}$. 

\begin{lemma}\label{pexpansion}
We have
\begin{equation*}
\begin{aligned}
\frac{1}{(2\pi\mathrm{i})^2}\,\wp(z;\,[\tau,\,1])&=
\frac{1}{12}+\frac{q_z}{(1-q_z)^2}\\
&+\sum_{n=1}^\infty
\left(
\frac{q_\tau^nq_z}{(1-q_\tau^nq_z)^2}+
\frac{q_\tau^nq_z^{-1}}{(1-q_\tau^nq_z^{-1})^2}\right)
-2\sum_{n=1}^\infty\frac{nq_\tau^n}{1-q_\tau^n}
\quad(z\in\mathbb{C},~\tau\in\mathbb{H}).
\end{aligned}
\end{equation*}
\end{lemma}
\begin{proof}
See \cite[Proposition 3 in Chapter 4]{Lang}. 
\end{proof}

The elliptic modular function $j$ on $\mathbb{H}$ is defined by
\begin{equation}\label{jdef}
j(\tau)=1728\,\frac{g_2(\tau)^3}{\Delta(\tau)}
=1728+46656\,\frac{g_3(\tau)^2}{\Delta(\tau)}
\quad(\tau\in\mathbb{H}),
\end{equation}
which is a weakly holomorphic modular function for $\mathrm{SL}_2(\mathbb{Z})$ 
(\cite[p. 39]{Lang}). 
Furthermore, its Fourier expansion with respect to $q_\tau$ is given by
\begin{equation}\label{jFourier}
j(\tau)
=q_\tau^{-1}+744+196884q_\tau+21493760q_\tau^2+\cdots\in q_\tau^{-1}\mathbb{Z}[[q_\tau]]
\end{equation}
with integer Fourier coefficients 
(see \cite[Theorem 11.8]{Cox}). 
\par
For $k\in\{1,\,2,\,3\}$ and $\mathbf{v}=
\begin{bmatrix}v_1&v_2\end{bmatrix}\in M_{1,\,2}(\mathbb{Q})\setminus
M_{1,\,2}(\mathbb{Z})$, the Fricke function $f^{(k)}_\mathbf{v}$
on $\mathbb{H}$ 
is defined by 
\begin{equation}\label{Frickedef}
f^{(k)}_\mathbf{v}(\tau)=
\left\{\begin{array}{ll}
\displaystyle\frac{g_2(\tau)g_3(\tau)}{\Delta(\tau)}\,
\wp(v_1\tau+v_2;\,[\tau,\,1]) & \textrm{if}~k=1,\vspace{0.1cm}\\
\displaystyle\frac{g_2(\tau)^2}{\Delta(\tau)}\,\wp(v_1\tau+v_2;\,[\tau,\,1])^2 & \textrm{if}~k=2,\vspace{0.1cm}\\
\displaystyle\frac{g_3(\tau)}{\Delta(\tau)}\,\wp(v_1\tau+v_2;\,[\tau,\,1])^3 & \textrm{if}~k=3
\end{array}\right.\qquad(\tau\in\mathbb{H}).
\end{equation}
Observe that $f^{(k)}_\mathbf{v}$ is holomorphic on $\mathbb{H}$. 
We simply write $f_\mathbf{v}$ for $f^{(1)}_\mathbf{v}$. By 
some basic properties of the Weierstrass $\wp$-function (see \cite[Lemma 10.4]{Cox}), we have
\begin{equation}\label{fufv}
f_\mathbf{u}=f_\mathbf{v}\quad\Longleftrightarrow\quad
\mathbf{u}\equiv\mathbf{v}~\textrm{or}~-\mathbf{v}\Mod{M_{1,\,2}(\mathbb{Z})}.
\end{equation}
For a positive integer $N$, consider the field $\mathcal{F}_N$ defined by
\begin{equation*}
\mathcal{F}_N=\left\{
\begin{array}{ll}
\mathbb{Q}(j)  & \textrm{if}~N=1,\\
\mathbb{Q}\left(j,\,f_\mathbf{v}\mid\mathbf{v}\in\displaystyle\frac{1}{N} M_{1,\,2}(\mathbb{Z})
\setminus M_{1,\,2}(\mathbb{Z})\right) & \textrm{if}~N\geq2.
\end{array}\right.
\end{equation*}

\begin{proposition}\label{modularGalois}
Let $N\geq2$ be an integer.
\begin{enumerate}
\item[\textup{(i)}]
The extension $\mathcal{F}_N/\mathcal{F}_1$ is Galois and
\begin{equation*}
\mathrm{Gal}(\mathcal{F}_N/\mathcal{F}_1)\simeq\mathrm{GL}_2(\mathbb{Z}/N\mathbb{Z})/\langle-I_2\rangle.
\end{equation*}
\item[\textup{(ii)}] More precisely, if $\alpha\in\mathrm{GL}_2(\mathbb{Z}/N\mathbb{Z})/\langle-I_2\rangle$, then 
\begin{equation}\label{fva}
f_\mathbf{v}^\alpha=f_{\mathbf{v}\widetilde{\alpha}}
\quad(\mathbf{v}\in\tfrac{1}{N}M_{1,\,2}(\mathbb{Z})\setminus
M_{1,\,2}(\mathbb{Z})),
\end{equation}
where $\widetilde{\alpha}\in M_2(\mathbb{Z})$ is any lift of $\alpha$. 
\item[\textup{(iii)}] Specifically, if $\gamma\in\mathrm{SL}_2(\mathbb{Z})$, then 
\begin{equation*}
f_\mathbf{v}\circ\gamma=f_\mathbf{v}^{[\gamma]}=f_{\mathbf{v}\gamma}
\quad(\mathbf{v}\in\frac{1}{N}M_{1,\,2}(\mathbb{Z})\setminus
M_{1,\,2}(\mathbb{Z})), 
\end{equation*}
where $[\gamma]$ is the reduction of $\gamma$ to $\mathrm{GL}_2(\mathbb{Z}/N\mathbb{Z})/
\langle-I_2\rangle$. 
\end{enumerate}
\end{proposition}
\begin{proof}
See  \cite[Theorem 6.6 (1), (2) and Proposition 6.21 (3)]{Shimura}. 
\end{proof}

We further note that $\mathcal{F}_N$ is the field
of meromorphic modular functions of level $N$ whose Fourier 
expansions with respect to $q_\tau^{1/N}$ have coefficients in the $N$th cyclotomic field
$\mathbb{Q}(\zeta_N)$; see \cite[Proposition 6.9 (1)]{Shimura}. 
\par
For $\mathbf{v}=
\begin{bmatrix}v_1&v_2\end{bmatrix}\in M_{1,\,2}(\mathbb{Q})\setminus
M_{1,\,2}(\mathbb{Z})$, the Siegel function $g_\mathbf{v}$
on $\mathbb{H}$ 
is given by the following infinite product formula
\begin{equation*}
\begin{aligned}
g_\mathbf{v}(\tau)=&-e^{\pi\mathrm{i}v_2(v_1-1)}
q_\tau^{(1/2)(v_1^2-v_1+1/6)}(1-q_\tau^{v_1}e^{2\pi\mathrm{i}v_2})\\
&\times\prod_{n=1}^\infty
(1-q_\tau^{n+v_1}e^{2\pi\mathrm{i}v_2})(1-q_\tau^{n-v_1}e^{-2\pi\mathrm{i}v_2})
\quad(\tau\in\mathbb{H}).
\end{aligned}
\end{equation*}
For the original definition of $g_\mathbf{v}$ in
terms of a Klein form and the Dedekind eta function, one may refer to
\cite[$\S$2.1]{K-L}. 
Note that $g_\mathbf{v}$ has no 
zeros or poles on $\mathbb{H}$. 
Furthermore, if $\mathbf{v}\in\frac{1}{N}M_{1,\,2}(\mathbb{Z})\setminus
M_{1,\,2}(\mathbb{Z})$, then $g_\mathbf{v}^{12N}$ belongs to $\mathcal{F}_N$ and 
satisfies the properties \eqref{fufv} and \eqref{fva}, with $f_\mathbf{v}$ replaced by 
$g_\mathbf{v}^{12N}$ (see \cite[Proposition 1.3 in Chapter 2]{K-L}). 
 
\begin{lemma}\label{difference}
If $\mathbf{u},\,\mathbf{v}\in M_{1,\,2}(\mathbb{Q})\setminus
M_{1,\,2}(\mathbb{Z})$ satisfy neither $\mathbf{u}\equiv\mathbf{v}\Mod{M_{1,\,2}(\mathbb{Z})}$ nor
$\mathbf{u}\equiv-\mathbf{v}\Mod{M_{1,\,2}(\mathbb{Z})}$, then 
\begin{equation*}
(f_\mathbf{u}-f_\mathbf{v})^6=\frac{j^2(j-1728)^3}{2^{30}3^{24}}\frac{g_{\mathbf{u}+\mathbf{v}}^6
g_{\mathbf{u}-\mathbf{v}}^6}
{g_\mathbf{u}^{12}g_\mathbf{v}^{12}}. 
\end{equation*}
\end{lemma}
\begin{proof}
See \cite[Theorem 2 in Chapter 18]{Lang} and \cite[p. 29 and p. 51]{K-L}. 
\end{proof}

\section {Primitive division points and Galois conjugates of their Weber values}

We describe primitive division points in terms of their annihilators
and relate the Galois conjugates of their Weber values to Fricke invariants. 
We then express these conjugates as values of Fricke functions at CM points
associated with reduced quadratic forms. 
\par
For $a\in\mathcal{O}_K\simeq\mathrm{End}(E)$ and $P=\phi_E(z+\mathcal{O}_K)\in E(\mathbb{C})$, 
we simply write $aP$ for $[a]_E(P)=\phi_E(az+\mathcal{O}_K)$. 

\begin{lemma}\label{primitivelemma}
Write $\mathfrak{m}=\prod_{\mathfrak{p}}
\mathfrak{p}^{e_\mathfrak{p}}$ ($e_\mathfrak{p}\geq0$), and 
let $P=\phi_E(z+\mathcal{O}_K)\in E(\mathbb{C})$ with $z\in K$.
\begin{enumerate}
\item[\textup{(i)}] We have
\begin{equation}\label{AnnP}
\mathrm{Ann}_{\mathcal{O}_K}(P)=\prod_{\mathfrak{p}}
\mathfrak{p}^{\max\{0,\,-v_\mathfrak{p}(z)\}},
\end{equation}
where $v_\mathfrak{p}$ is the normalized $\mathfrak{p}$-adic
valuation on $K$ with $v_\mathfrak{p}(0)=\infty$. Thus we obtain
\begin{equation*}
\textrm{$P$ is a primitive $\mathfrak{m}$-division point}
\quad\Longleftrightarrow\quad\max\{0,\,-v_\mathfrak{p}(z)\}=e_\mathfrak{p}~\textrm{for all}~\mathfrak{p}. 
\end{equation*} 
\item[\textup{(ii)}] Let $a\in\mathcal{O}_K$. If $P$ is a primitive $\mathfrak{m}$-division point, then 
\begin{equation*}
\textrm{$ aP$ is a primitive $\mathfrak{m}$-division point}
\quad\Longleftrightarrow\quad a+\mathfrak{m}\in(\mathcal{O}_K/\mathfrak{m})^\times. 
\end{equation*}
\end{enumerate}
\end{lemma}
\begin{proof}
\begin{enumerate}
\item[(i)] We deduce that for $a\in\mathcal{O}_K$
\begin{eqnarray*}
a\in\prod_{\mathfrak{p}}
\mathfrak{p}^{\max\{0,\,-v_\mathfrak{p}(z)\}}
&\Longleftrightarrow&
v_\mathfrak{p}(a)
\geq\max\{0,\,-v_\mathfrak{p}(z)\}~\textrm{for all}~\mathfrak{p}\\
&\Longleftrightarrow&
a\in\mathcal{O}_K~\textrm{and}~az\in\mathcal{O}_K\\
&\Longleftrightarrow& aP=O_E\\
&\Longleftrightarrow&a\in\mathrm{Ann}_{\mathcal{O}_K}(P).  
\end{eqnarray*} 
This proves \eqref{AnnP}. The second assertion follows from the definition of a 
primitive $\mathfrak{m}$-division point. 
\item[(ii)]
Consider the multiplication map
\begin{equation*}
\mu_a:\mathcal{O}_K/\mathfrak{m}\rightarrow\mathcal{O}_K/\mathfrak{m},\quad
b+\mathfrak{m}\mapsto ab+\mathfrak{m},
\end{equation*}
which is an $\mathcal{O}_K$-module homomorphism. We then
derive by
the assumption $\mathrm{Ann}_{\mathcal{O}_K}(P)=\mathfrak{m}$
that for $b\in\mathcal{O}_K$
\begin{align*}
b+\mathfrak{m}\in\mathrm{ker}\,\mu_a
~\Leftrightarrow~
ab\in\mathfrak{m}
~\Leftrightarrow~(ab)P=O_E
~\Leftrightarrow~b( aP)=O_E
~\Leftrightarrow~b\in\mathrm{Ann}_{\mathcal{O}_K}(aP),
\end{align*}
and hence
\begin{equation}\label{kernel}
\mathrm{Ann}_{\mathcal{O}_K}( aP)\supseteq\mathfrak{m}\quad
\textrm{and}\quad\mathrm{ker}\,\mu_a=\mathrm{Ann}_{\mathcal{O}_K}( aP)/\mathfrak{m}. 
\end{equation}
Therefore we achieve that
\begin{align*}
\textrm{$ aP$ is a primitive $\mathfrak{m}$-division point}&\Longleftrightarrow
\mathrm{Ann}_{\mathcal{O}_K}( aP)=\mathfrak{m}\\
&\Longleftrightarrow\mathrm{ker}\,\mu_a=\{0+\mathfrak{m}\}\quad\textrm{by \eqref{kernel}}\\
&\Longleftrightarrow\textrm{$\mu_a$ is injective}\\
&\Longleftrightarrow\textrm{$\mu_a$ is bijective}\quad\textrm{since}~|\mathcal{O}_K/\mathfrak{m}|
=\mathrm{N}(\mathfrak{m})<\infty\\
&\Longleftrightarrow\textrm{$a+\mathfrak{m}$ is a unit in $\mathcal{O}_K/\mathfrak{m}$}. 
\end{align*}
\end{enumerate}
\end{proof}

\begin{remark}
 Lemma \ref{primitivelemma} (ii) also follows directly from (i). 
\end{remark}

\begin{definition}\label{invariants}
Let $N_\mathfrak{m}>1$ be the least positive integer in $\mathfrak{m}$. 
For $C\in C(\mathfrak{m})$, take any ideal $\mathfrak{c}$ of $\mathcal{O}_K$ in $C$, and express
\begin{align*}
\mathfrak{m}\mathfrak{c}^{-1}&=[\omega_1,\,\omega_2]\quad
\textrm{with}~\omega_1,\,\omega_2\in\mathbb{C}~\textrm{such that}~\omega=\frac{\omega_1}{\omega_2}\in
\mathbb{H},\\
1&=r_1\omega_1+r_2\omega_2\quad\textrm{with}~r_1,\,r_2\in\frac{1}{N_\mathfrak{m}}\mathbb{Z}. 
\end{align*}
We define the \textit{Fricke invariant} 
$f_\mathfrak{m}(C)$ and the \textit{Siegel-Ramachandra invariant} $g_\mathfrak{m}(C)$ by 
\begin{equation*}
f_\mathfrak{m}(C)=f_{\left[\begin{smallmatrix}r_1& r_2\end{smallmatrix}\right]}^{(|\mathcal{O}_K^\times|/2)}
(\omega)\quad\textrm{and}\quad
g_\mathfrak{m}(C)=g_{\left[\begin{smallmatrix}r_1& r_2\end{smallmatrix}\right]}(\omega)^{12N_\mathfrak{m}},
\end{equation*}
respectively. These values depend only on the ray class $C$, not on
the choices of $\mathfrak{c}$, $\omega_1$ and $\omega_2$
(see \cite[$\S$2.1 and $\S$11.1]{K-L} and \cite[$\S$6.2 and $\S$6.3]{Lang}). 
\end{definition}

Let $\sigma_\mathfrak{m}:C(\mathfrak{m})\rightarrow\mathrm{Gal}(K_\mathfrak{m}/K)$
denote the Artin map modulo $\mathfrak{m}$. 
By the Shimura reciprocity law, 
we have the following result. 

\begin{proposition}\label{fginvariants}
With the notation of \textup{Definition \ref{invariants}}, 
the invariants $f_\mathfrak{m}(C)$ and $g_\mathfrak{m}(C)$ belong to $K_\mathfrak{m}$. 
For every $C'\in C(\mathfrak{m})$, they satisfy 
\begin{equation*}
f_\mathfrak{m}(C)^{\sigma_\mathfrak{m}(C')}=f_\mathfrak{m}(CC')\quad\textrm{and}\quad
g_\mathfrak{m}(C)^{\sigma_\mathfrak{m}(C')}=g_\mathfrak{m}(CC').
\end{equation*}
\end{proposition}
\begin{proof}
See \cite[Theorem 1.1 in Chapter 11]{K-L}. 
\end{proof}

For a fractional ideal $\mathfrak{a}$ relatively prime to $\mathfrak{m}$, let $[\mathfrak{a}]_\mathfrak{m}$ denote its class in $C(\mathfrak{m})$. 

\begin{lemma}\label{hconjugate}
Let $P=\phi_E(z+\mathcal{O}_K)\in E$ with $z\in K$ be a primitive $\mathfrak{m}$-division point,
and let $a\in\mathcal{O}_K$.
If $a+\mathfrak{m}\in(\mathcal{O}_K/\mathfrak{m})^\times$, then 
\begin{equation*}
h_E(P)^{\sigma_\mathfrak{m}([a\mathcal{O}_K]_\mathfrak{m})}=h_E( aP). 
\end{equation*}
\end{lemma}
\begin{proof}
Note that $\mathfrak{c}=z\mathfrak{m}$ is a nonzero ideal of $\mathcal{O}_K$
which is relatively prime to $\mathfrak{m}$ by the assumption
$\mathrm{Ann}_{\mathcal{O}_K}(P)=\mathfrak{m}$ and Lemma \ref{primitivelemma} (i). 
Set 
\begin{equation*}
C=[\mathfrak{c}]_\mathfrak{m}=[z\mathfrak{m}]_\mathfrak{m}\quad\textrm{and}\quad D=[az\mathfrak{m}]_\mathfrak{m}. 
\end{equation*}
Then we see by \eqref{Weberh}, \eqref{Frickedef} and Definition \ref{invariants} that
\begin{align}
h_E(P)&=f_\mathfrak{m}(C),\label{hf1}\\
h_E( aP)&=f_\mathfrak{m}(D).\label{hf2}
\end{align}
And, we find that
\begin{align*}
h_E(P)^{\sigma_\mathfrak{m}([a\mathcal{O}_K]_\mathfrak{m})}&=f_\mathfrak{m}([z\mathfrak{m}]_\mathfrak{m})^{\sigma_\mathfrak{m}([a\mathcal{O}_K]_\mathfrak{m})}\quad\textrm{by \eqref{hf1}}\\
&=f_\mathfrak{m}([az\mathfrak{m}]_\mathfrak{m})\quad\textrm{by Proposition \ref{fginvariants}}\\
&=h_E( aP)\quad\textrm{by \eqref{hf2}}. 
\end{align*}
\end{proof}

\begin{remark}\label{CP}
Let $P=\phi_E(z+\mathcal{O}_K)\in E(\mathbb{C})$ with $z\in K$ be a primitive $\mathfrak{m}$-division point.
\begin{enumerate}
\item[(i)] In the proof of Lemma \ref{hconjugate}, one can readily check that the ray class $C=[z\mathfrak{m}]_\mathfrak{m}$ depends only 
on $P$, not on the choice of $z$ (cf. \cite[$\S$V.1]{Janusz}). So we can write $C_P$ for $C$. 
By \eqref{O_K}, \eqref{hf1} and Definition \ref{invariants}, one sees that
\begin{equation}\label{fCP}
h_E(P)=f_\mathfrak{m}(C_P)=f_\mathbf{v}^{(|\mathcal{O}_K^\times|/2)}(\tau_K)\quad\textrm{for}~
\mathbf{v}=\begin{bmatrix}v_1&v_2
\end{bmatrix}\in M_{1,\,2}(\mathbb{Q})~\textrm{satisfying}~
z=v_1\tau_K+v_2,
\end{equation}
and
\begin{equation}\label{gCP}
g_\mathfrak{m}(C_P)=g_\mathbf{v}(\tau_K)^{12N_\mathfrak{m}}.
\end{equation}
\item[(ii)] Moreover, if $n$ is an integer such that $n+\mathfrak{m}\in(\mathcal{O}_K/\mathfrak{m})^\times$, then 
we get by Definition \ref{invariants} and Proposition \ref{fginvariants} that
\begin{align*}
f_\mathfrak{m}(C_P)^{\sigma_\mathfrak{m}([n\mathcal{O}_K]_\mathfrak{m})}
=f_\mathfrak{m}([z\mathfrak{m}]_\mathfrak{m})^{\sigma_\mathfrak{m}([n\mathcal{O}_K]_\mathfrak{m})}
=f_\mathfrak{m}([nz\mathfrak{m}]_\mathfrak{m})=f_{n\mathbf{v}}^{(|\mathcal{O}_K^\times|/2)}(\tau_K).
\end{align*}
In a similar way, we establish that
\begin{equation}\label{gsg}
g_\mathfrak{m}(C_P)^{\sigma_\mathfrak{m}([n\mathcal{O}_K]_\mathfrak{m})}
=g_{n\mathbf{v}}(\tau_K)^{12N_\mathfrak{m}}\quad(n\in\mathbb{Z}~\textrm{satisfying}~n+\mathfrak{m}
\in(\mathcal{O}_K/\mathfrak{m})^\times). 
\end{equation}
\item[(iii)] Lemma \ref{hconjugate}
also follows from the idelic argument in the proof of \cite[Theorem 7 in Chapter 10]{Lang}. 
\end{enumerate}
\end{remark}

Let $\mathcal{Q}(d_K)$ be the 
set of primitive positive definite binary quadratic forms over $\mathbb{Z}$ of discriminant
$d_K$, that is,
\begin{equation*}
\mathcal{Q}(d_K)=\{ax^2+bxy+cy^2\in\mathbb{Z}[x,\,y]\mid
\gcd(a,\,b,\,c)=1,~b^2-4ac=d_K,\,a>0\}.
\end{equation*}
For each $Q=ax^2+bxy+cy^2\in\mathcal{Q}(d_K)$,
define $\omega_Q$ as the zero of the polynomial $Q(x,\,1)$ lying in $\mathbb{H}$, 
namely,
\begin{equation*}
\omega_Q=\frac{-b+\sqrt{d_K}}{2a}. 
\end{equation*}
We say that a form $Q=ax^2+bxy+cy^2\in\mathcal{Q}(d_K)$ is reduced if
\begin{equation}\label{reducedcondition}
|b|\leq a\leq c\quad\textrm{and}\quad
b\geq0~\textrm{if either}~|b|=a~\textrm{or}~a=c. 
\end{equation}
In this case, one sees directly from \eqref{reducedcondition} that
\begin{equation}\label{abound}
a\leq\sqrt{\frac{|d_K|}{3}}.
\end{equation}
Moreover, we call $Q$ the principal reduced form if $a=1$, in which case
\begin{equation*}
Q=\left\{\begin{array}{ll}
\displaystyle x^2-\frac{d_K}{4}y^2 & \textrm{if}~d_K\equiv0\Mod{4},\\
\displaystyle x^2+xy+\frac{1-d_K}{4}y^2 & \textrm{if}~d_K\equiv1\Mod{4}. 
\end{array}\right.
\end{equation*} 
Let $\mathcal{Q}_\mathrm{red}(d_K)$ denote
the set of reduced forms in $\mathcal{Q}(d_K)$. 
As is well known, 
\begin{equation}\label{bijective}
\mathcal{Q}_\mathrm{red}(d_K)\rightarrow C(\mathcal{O}_K)
~\textrm{sending}~
Q~\textrm{to}~[[\omega_Q,\,1]]_{\mathcal{O}_K}~
\textrm{is bijective}
\end{equation} 
(see \cite[Theorems 2.8 and 7.7]{Cox}).

\begin{lemma}\label{fenough}
With the notation of \textup{Definition \ref{invariants}}, we have
\begin{equation*}
f_\mathfrak{m}(C)=f_\mathbf{w}^{(|\mathcal{O}_K^\times|/2)}(\omega_Q)\quad\textrm{and}\quad
g_\mathfrak{m}(C)=g_{\mathbf{w}}(\omega_Q)^{12N_\mathfrak{m}}
\end{equation*}
for some $Q\in\mathcal{Q}_\mathrm{red}(d_K)$
and $\mathbf{w}\in\frac{1}{N_\mathfrak{m}}M_{1,\,2}(\mathbb{Z})
\setminus M_{1,\,2}(\mathbb{Z})$.
\end{lemma}
\begin{proof}
We get by \eqref{bijective} that
\begin{equation*}
\mathfrak{m}\mathfrak{c}^{-1}=\nu\,[\omega_Q,\,1]\quad\textrm{for some}~
\nu\in K^\times~\textrm{and}~
Q\in\mathcal{Q}_\mathrm{red}(d_K).
\end{equation*}
Hence we attain by Definition \ref{invariants} that
\begin{equation*}
f_\mathfrak{m}(C)=f_{\mathbf{w}}^{(|\mathcal{O}_K^\times|/2)}(\omega_Q)\quad\textrm{and}\quad
g_\mathfrak{m}(C)=g_{\mathbf{w}}(\omega_Q)^{12N_\mathfrak{m}}
\end{equation*}
for some $Q\in\mathcal{Q}_\mathrm{red}(d_K)$
and $\mathbf{w}\in\frac{1}{N_\mathfrak{m}}M_{1,\,2}(\mathbb{Z})
\setminus M_{1,\,2}(\mathbb{Z})$.
\end{proof}

\section {The second Kronecker limit formula}

We first recall the second Kronecker limit formula,
which relates character sums
of logarithms of Siegel-Ramachandra invariants to 
Hecke $L$-values at $s=1$. We then  
derive an identity for differences of Weber values used in the 
character argument of Section \ref{sect6}. 
\par
Let $\chi$ be a nonprincipal character of the ray class group $C(\mathfrak{m})$.
We define the Stickelberger element $S_\mathfrak{m}(\chi)$ 
by
\begin{equation}\label{Sdef}
S_\mathfrak{m}(\chi)=\sum_{C\in C(\mathfrak{m})}\chi(C)
\ln|g_\mathfrak{m}(C)|.
\end{equation}
Furthermore, we define the $L$-function $L_\mathfrak{m}(\cdot,\,\chi)$ by 
\begin{equation*}
L_\mathfrak{m}(s,\,\chi)=\sum_{\mathfrak{a}}\frac{\chi([\mathfrak{a}]_\mathfrak{m})}
{\mathrm{N}(\mathfrak{a})^s}\quad(s\in\mathbb{C}~\textrm{with}~\mathrm{Re}(s)>1),
\end{equation*}
where $\mathfrak{a}$ runs over all nonzero ideals of $\mathcal{O}_K$ relatively prime to $\mathfrak{m}$.
Then it is well known that
$L_\mathfrak{m}(\cdot,\,\chi)$ extends holomorphically to
a neighborhood of $s=1$, and 
\begin{equation}\label{Lnonzero}
L_\mathfrak{m}(1,\,\chi)\neq0 
\end{equation}
(see \cite[Theorem 2.14, Proposition 4.1 in Chapter IV and
Theorem 10.2 in Chapter V]{Janusz}). 
\par
Denote by $\mathfrak{f}_\chi$ the conductor of $\chi$, namely,
the smallest divisor of $\mathfrak{m}$
such that $\chi$ factors through the natural projection
$C(\mathfrak{m})\twoheadrightarrow C(\mathfrak{f}_\chi)$. 

\begin{proposition}\label{LS}
Let $\chi_0$ be the primitive character on $C(\mathfrak{f}_\chi)$ associated with $\chi$. 
If $\mathfrak{f}_\chi\neq\mathcal{O}_K$, then we have
\begin{equation*}
L_{\mathfrak{f}_\chi}(1,\,\chi_0)\left(\prod_{\mathfrak{p}\mid\mathfrak{m},~\mathfrak{p}\nmid\mathfrak{f}_\chi}(1-\overline{\chi}_0([\mathfrak{p}]))\right)=
-\frac{\pi\chi_0([\gamma\sqrt{d_K}\mathfrak{f}_\chi])}
{3N_\mathfrak{m}\sqrt{|d_K|}\omega(\mathfrak{m})T_\gamma
(\overline{\chi}_0)}\,S_\mathfrak{m}(\overline{\chi}),
\end{equation*}
where $\gamma$ is an element of $K$ so that $\gamma\sqrt{d_K}\mathfrak{f}_\chi$ is 
a nonzero ideal of $\mathcal{O}_K$ relatively prime to $\mathfrak{f}_\chi$,
$N_\mathfrak{m}$ is the least positive integer in $\mathfrak{m}$,
\begin{align*}
\omega(\mathfrak{m})=|\{\nu\in\mathcal{O}_K^\times\mid\nu\equiv 1\Mod{\mathfrak{m}}\}|,
\quad
T_\gamma(\overline{\chi}_0)=\sum_{\mu+\mathfrak{f}_\chi\in
(\mathcal{O}_K/\mathfrak{f}_\chi)^\times}
\overline{\chi}_0([\mu\mathcal{O}_K])
e^{2\pi\mathrm{i}\mathrm{Tr}_{K/\mathbb{Q}}(\mu\gamma)}. 
\end{align*}
\end{proposition}
\begin{proof}
See \cite[Theorem 2.1 in Chapter 11]{K-L} or \cite[Theorem 9 in Chapter II]{Siegel}. 
\end{proof}

\begin{remark}\label{barremark}
By the definition \eqref{Sdef}, we see that $S_\mathfrak{m}(\overline{\chi})=
\overline{S_\mathfrak{m}(\chi)}$. 
\end{remark}

\begin{lemma}\label{Schilemma}
Let $\mathfrak{m}$ be a nonzero proper ideal of $\mathcal{O}_K$ which is one of the forms
\begin{equation}\label{4forms}
\mathfrak{p}_{3,\,\mathrm{split}}^2,\quad
\mathfrak{p}_{3,\,\mathrm{split}}^2\mathfrak{p}'_{3,\,\mathrm{split}},\quad
\mathfrak{p}_{3,\,\mathrm{noninert}}\mathfrak{p}_{5,\,\mathrm{noninert}},\quad
\mathfrak{p}_{3,\,\mathrm{split}}\mathfrak{p}'_{3,\,\mathrm{split}}\mathfrak{p}_{5,\,\mathrm{noninert}},
\end{equation}
so $2+\mathfrak{m}\in(\mathcal{O}_K/\mathfrak{m})^\times$. 
Let $P\in E(\mathbb{C})$ be a primitive $\mathfrak{m}$-division point, and let
\begin{equation*}
\delta=h_E(2P)-h_E(P).
\end{equation*}
Let $\chi$ be a character of $C(\mathfrak{m})$.
If $h_K>1$ and
$\mathfrak{f}_\chi=\mathfrak{m}$, then we have $\delta\neq0$ and
\begin{equation*}
\sum_{C\in C(\mathfrak{m})}\chi(C)\log|(\delta^{12N_\mathfrak{m}})^{\sigma_\mathfrak{m}(C)}|
=-\overline{\chi}(C_P)(1+2\overline{\chi}([2\mathcal{O}_K]_\mathfrak{m}))S_\mathfrak{m}(\chi), 
\end{equation*}
where $C_P$ is the ray class described in \textup{Remark \ref{CP} (i)}. 
\end{lemma}
\begin{proof}
Since $h_K>1$, we get $|\mathcal{O}_K^\times|/2=1$ and
\begin{equation}\label{not01728}
j(\tau_K)\neq0,\,1728\quad(\textrm{if}~h_K>1)
\end{equation}
(see \cite[p. 167, Theorems 11.2 and 12.34]{Cox}). By 
\eqref{fCP}, we achieve
\begin{equation}\label{hPfv}
h_E(P)=f_\mathbf{v}(\tau_K)\quad\textrm{for some}~\mathbf{v}\in M_{1,\,2}(\mathbb{Q})\setminus
M_{1,\,2}(\mathbb{Z}).
\end{equation}
Here $\mathbf{v}=\begin{bmatrix}v_1&v_2\end{bmatrix}$
is chosen so that $P=\phi_E(z+\mathcal{O}_K)$, where $z=v_1\tau_K+v_2$. 
Note that $3\mathbf{v}\not\in M_{1,\,2}(\mathbb{Z})$
because $\mathfrak{m}$ is one of the ideals listed in \eqref{4forms}. 
Hence we find that
\begin{equation}\label{d6}
\begin{aligned}
\delta^6&=\left(f_{2\mathbf{v}}(\tau_K)-f_\mathbf{v}(\tau_K)\right)^6\quad\textrm{by \eqref{hPfv}}\\
&=J_K\frac{g_{3\mathbf{v}}(\tau_K)^6}
{g_{2\mathbf{v}}(\tau_K)^{12}g_\mathbf{v}(\tau_K)^6}\quad\textrm{by Lemma \ref{difference},}\\
&\hspace{4.5cm}\textrm{where}~
J_K=\frac{j(\tau_K)^2(j(\tau_K)-1728)^3}{2^{30}3^{24}}\in H_K\\
&\neq0\quad\textrm{by \eqref{not01728}}.
\end{aligned}
\end{equation}
Let $\mathfrak{m}'=\mathrm{Ann}_{\mathcal{O}_K}(3P)$. 
Then we have by Lemma \ref{primitivelemma} (i)  that
\begin{equation*}
\mathfrak{m}'=\frac{\mathfrak{m}}{
\gcd(\mathfrak{m},\,3\mathcal{O}_K)},
\end{equation*}
and so 
\begin{equation*}
\mathfrak{m}'\mid\mathfrak{m},\quad
\mathfrak{m}'\neq\mathfrak{m},\quad\mathfrak{m}'\neq\mathcal{O}_K. 
\end{equation*}
Observe that 
\begin{equation}\label{vanish}
\begin{aligned}
&\sum_{C\in C(\mathfrak{m})}\chi(C)\log|\nu^{\sigma_\mathfrak{m}(C)}|\quad\textrm{for any}~\nu\in K_{\mathfrak{m}'}^\times\\
=&\sum_{A\in\mathrm{ker}(C(\mathfrak{m})
\rightarrow C(\mathfrak{m}'))}
\sum_{B\in \mathcal{R}}
\chi(AB)\log|\nu^{\sigma_\mathfrak{m}(AB)}|\\
&\hspace{1cm}
\textrm{where $\mathcal{R}$ is a complete set of representatives
of}~
C(\mathfrak{m})/\mathrm{ker}(C(\mathfrak{m})
\rightarrow C(\mathfrak{m}'))\\
=&\left(\sum_{A\in\mathrm{ker}(C(\mathfrak{m})
\rightarrow C(\mathfrak{m}'))}\chi(A)\right)\cdot
\left(\sum_{B\in \mathcal{R}}
\chi(B)\log|\nu^{\sigma_\mathfrak{m}(B)}|\right)
\\
&\hspace{8cm}
\textrm{since}~\nu\in K_{\mathfrak{m}'}~\textrm{and so}~
\nu^{\sigma_\mathfrak{m}(A)}=\nu\\
=&0\cdot\left(\sum_{B\in \mathcal{R}}
\chi(B)\log|\nu^{\sigma_\mathfrak{m}(B)}|\right)\quad\textrm{by \cite[Proposition 4 in Chapter VI]{Serre},}
\\&\hspace{3cm}
\textrm{because}~\mathfrak{f}_\chi=\mathfrak{m}~\textrm{implies that
$\chi$ is nontrivial on $\mathrm{ker}(C(\mathfrak{m})
\rightarrow C(\mathfrak{m}'))$}\\
=&0.
\end{aligned}
\end{equation}
Now, we derive that
\begin{align*}
&\sum_{C\in C(\mathfrak{m})}\chi(C)\log|(\delta^{12N_\mathfrak{m}})^{\sigma_\mathfrak{m}(C)}|\\
=&~2N_\mathfrak{m}\sum_{C\in C(\mathfrak{m})}\chi(C)\log|J_K^{\sigma_\mathfrak{m}(C)}|
+\sum_{C\in C(\mathfrak{m})}\chi(C)\log|(g_{3\mathbf{v}}(\tau_K)^{12N_\mathfrak{m}})^{\sigma_\mathfrak{m}(C)}|\\
&-2\sum_{C\in C(\mathfrak{m})}\chi(C)\log|(g_{2\mathbf{v}}(\tau_K)^{12N_\mathfrak{m}})^{\sigma_\mathfrak{m}(C)}|
-\sum_{C\in C(\mathfrak{m})}\chi(C)\log|(g_\mathbf{v}(\tau_K)^{12N_\mathfrak{m}})^{\sigma_\mathfrak{m}(C)}|
\quad\textrm{by \eqref{d6}}\\\
=&~2N_\mathfrak{m}\sum_{C\in C(\mathfrak{m})}\chi(C)\log|J_K^{\sigma_\mathfrak{m}(C)}|
+\sum_{C\in C(\mathfrak{m})}\chi(C)\log|(g_{\mathfrak{m}'}([3z\mathfrak{m}']_{\mathfrak{m}'})^{N_\mathfrak{m}/N_{\mathfrak{m}'}})^{\sigma_\mathfrak{m}(C)}|\\
&-2\sum_{C\in C(\mathfrak{m})}\chi(C)\log|(g_\mathfrak{m}(C_P)^
{\sigma_\mathfrak{m}
([2\mathcal{O}_K]_\mathfrak{m})})^{\sigma_\mathfrak{m}(C)}|
-\sum_{C\in C(\mathfrak{m})}\chi(C)\log|g_\mathfrak{m}(C_P)^{\sigma_\mathfrak{m}(C)}|\\
&\hspace{12cm}\textrm{by \eqref{gCP} and \eqref{gsg}}\\
=&0+0-2\overline{\chi}(C_P[2\mathcal{O}_K]_\mathfrak{m})S_\mathfrak{m}(\chi)
-\overline{\chi}(C_P)S_\mathfrak{m}(\chi)\\
&\hspace{7cm}\textrm{by \eqref{vanish} due to the fact}~J_K,\,g_{\mathfrak{m}'}([3z\mathfrak{m}']_{\mathfrak{m}'})
\in K_{\mathfrak{m}'}^\times\\
=&-\overline{\chi}(C_P)(1+2\overline{\chi}([2\mathcal{O}_K]_\mathfrak{m}))S_\mathfrak{m}(\chi).
\end{align*}

\end{proof}

\section {Proof using an $L$-function argument}\label{sect6}

In this section, we prove Theorem \ref{main} for four
of the seven families in Proposition \ref{residual}
using the second Kronecker limit formula. 

\begin{lemma}\label{charlemma}
Let $\mathfrak{m}$ be a nonzero proper ideal of $\mathcal{O}_K$ 
which is one of the forms 
\begin{equation*}
\mathfrak{p}_{3,\,\mathrm{split}}^2,\quad
\mathfrak{p}_{3,\,\mathrm{split}}^2\mathfrak{p}'_{3,\,\mathrm{split}},\quad
\mathfrak{p}_{3,\,\mathrm{noninert}}\mathfrak{p}_{5,\,\mathrm{noninert}},\quad
\mathfrak{p}_{3,\,\mathrm{split}}\mathfrak{p}'_{3,\,\mathrm{split}}\mathfrak{p}_{5,\,\mathrm{noninert}}.
\end{equation*}
 Let $P\in E(\mathbb{C})$ be a primitive $\mathfrak{m}$-division point, 
and let $L=K(h_E(P))$.
If $L\neq K_\mathfrak{m}$,
then there exists a character $\chi$ of $C(\mathfrak{m})$
such that
\begin{equation}\label{required}
\mathfrak{f}_\chi=\mathfrak{m}\quad\textrm{and}\quad
\chi|_{\sigma_\mathfrak{m}^{-1}(\mathrm{Gal}(K_\mathfrak{m}/L))}\neq1. 
\end{equation}
\end{lemma}
\begin{proof}
If $h_K=1$, then \eqref{overH} gives $L=K_\mathfrak{m}$. So we may assume that $h_K>1$, and hence
$\mathcal{O}_K^\times=\{\pm1\}$. 
The assumption $L\neq K_\mathfrak{m}$ yields
\begin{equation}\label{nontrivialL}
\sigma_\mathfrak{m}^{-1}(\mathrm{Gal}(K_\mathfrak{m}/L))\neq\{[\mathcal{O}_K]_\mathfrak{m}\}.
\end{equation}
Let
\begin{equation*}
R_\mathfrak{m}:=\mathrm{ker}(C(\mathfrak{m})\twoheadrightarrow C(\mathcal{O}_K)). 
\end{equation*}
Since
$R_\mathfrak{m}
\simeq\mathrm{Gal}(K_\mathfrak{m}/H_K)$ and 
$H_KL=K_\mathfrak{m}$ by \eqref{overH},
we observe that
\begin{equation}\label{intersection}
R_\mathfrak{m}
\cap\sigma_\mathfrak{m}^{-1}(\mathrm{Gal}(K_\mathfrak{m}/L))=\{[\mathcal{O}_K]_\mathfrak{m}\}.
\end{equation}
\par
On the other hand, we have the natural isomorphisms
\begin{equation}\label{ZZOP}
\mathbb{Z}/3\mathbb{Z}\simeq\mathcal{O}_K/\mathfrak{p}_{3,\,\mathrm{noninert}},\quad
\mathbb{Z}/9\mathbb{Z}\simeq\mathcal{O}_K/\mathfrak{p}_{3,\,\mathrm{split}}^2,
\quad
\mathbb{Z}/5\mathbb{Z}\simeq\mathcal{O}_K/\mathfrak{p}_{5,\,\mathrm{noninert}}. 
\end{equation}
Note that the multiplicative groups $(\mathbb{Z}/3\mathbb{Z})^\times$,
$(\mathbb{Z}/9\mathbb{Z})^\times$, $(\mathbb{Z}/5\mathbb{Z})^\times$
are cyclic groups of order $2$, $6$, $4$
generated by $-1+3\mathbb{Z}$, $2+9\mathbb{Z}$, $2+5\mathbb{Z}$, respectively. 
Thus we define characters
$\xi_1$, $\xi_2$, $\xi_3$ on $(\mathbb{Z}/3\mathbb{Z})^\times$,
$(\mathbb{Z}/9\mathbb{Z})^\times$, $(\mathbb{Z}/5\mathbb{Z})^\times$ by letting
\begin{equation}\label{x1x2x3}
\xi_1(-1+3\mathbb{Z})=-1,\quad \xi_2(2+9\mathbb{Z})=\zeta_6,\quad
\xi_3(2+5\mathbb{Z})=\mathrm{i},
\end{equation}
respectively. We then see that
\begin{equation}\label{-1}
\xi_1(-1+3\mathbb{Z})=\xi_2(-1+9\mathbb{Z})=\xi_3(-1+5\mathbb{Z})=-1. 
\end{equation}
If we write $\mathfrak{m}=\prod_{\mathfrak{p}\mid\mathfrak{m}}\mathfrak{p}^{e_\mathfrak{p}}$,
then we get by the Chinese remainder theorem that
\begin{equation*}
(\mathcal{O}_K/\mathfrak{m})^\times\simeq\prod_{\mathfrak{p}\mid\mathfrak{m}}
(\mathcal{O}_K/\mathfrak{p}^{e_\mathfrak{p}})^\times.
\end{equation*}
Thus, by \eqref{ZZOP},  we can identify $(\mathcal{O}_K/\mathfrak{m})^\times$ 
with the group in the second column 
of Table \ref{table1},
and define a character $\theta$ as follows: 
\begin{center}
\begin{tabular}{@{}ccc@{}}
\toprule
$\mathfrak{m}$ & Identification of $(\mathcal{O}_K/\mathfrak{m})^\times$ 
  & Definition of $\theta$ \\
\midrule\midrule
$\mathfrak{p}_{3,\,\mathrm{split}}^2$ & $(\mathbb{Z}/9\mathbb{Z})^\times$
  & $\theta(x)=\xi_2(x)^2$ \\
$\mathfrak{p}_{3,\,\mathrm{split}}^2\mathfrak{p}'_{3,\,\mathrm{split}}$ & $(\mathbb{Z}/9\mathbb{Z})^\times\times(\mathbb{Z}/3\mathbb{Z})^\times$
  & $\theta(x,\,y)=\xi_2(x)\xi_1(y)$ \\
$\mathfrak{p}_{3,\,\mathrm{noninert}}\mathfrak{p}_{5,\,\mathrm{noninert}}$ & 
$(\mathbb{Z}/3\mathbb{Z})^\times\times(\mathbb{Z}/5\mathbb{Z})^\times$
  & $\theta(x,\,y)=\xi_1(x)\xi_3(y)$ \\
$\mathfrak{p}_{3,\,\mathrm{split}}\mathfrak{p}'_{3,\,\mathrm{split}}\mathfrak{p}_{5,\,\mathrm{noninert}}$ & 
$(\mathbb{Z}/3\mathbb{Z})^\times\times
(\mathbb{Z}/3\mathbb{Z})^\times\times(\mathbb{Z}/5\mathbb{Z})^\times$
  & $\theta(x,\,y,\,z)=\xi_1(x)\xi_1(y)\xi_3(z)^2$ \\
\bottomrule
\end{tabular}
\captionof{table}{A character $\theta$ of $(\mathcal{O}_K/\mathfrak{m})^\times$}\label{table1}
\end{center}
The exact sequence given in \eqref{exact} yields the surjective homomorphism
\begin{equation*}
q_\mathfrak{m}:(\mathcal{O}_K/\mathfrak{m})^\times\twoheadrightarrow R_\mathfrak{m}\simeq\mathrm{Gal}(K_\mathfrak{m}/H_K),
\quad a+\mathfrak{m}\mapsto [a\mathcal{O}_K]_\mathfrak{m}
\end{equation*}
with
\begin{equation*}
\quad\textrm{ker}(q_\mathfrak{m})=
\textrm{the image of $\mathcal{O}_K^\times$ in 
$(\mathcal{O}_K/\mathfrak{m})^\times$}=\{1+\mathfrak{m},\,-1+\mathfrak{m}\}. 
\end{equation*}
Since the character $\theta$ of $(\mathcal{O}_K/\mathfrak{m})^\times$
is trivial on $\mathrm{ker}(q_\mathfrak{m})=\{1+\mathfrak{m},\,
-1+\mathfrak{m}\}$ by \eqref{-1}, it induces a character
$\widetilde{\theta}$ of $(\mathcal{O}_K/\mathfrak{m})^\times/\mathrm{ker}(q_\mathfrak{m})$.
Moreover, since
$(\mathcal{O}_K/\mathfrak{m})^\times/\mathrm{ker}(q_\mathfrak{m})$ 
injects into the finite abelian group $C(\mathfrak{m})$ by \eqref{exact}, 
$\widetilde{\theta}$ extends to a character $\eta$ of $C(\mathfrak{m})$ (see
\cite[Proposition 1 in Chapter VI]{Serre}). 
\par
For a divisor $\mathfrak{a}$ of $\mathfrak{m}$, let
\begin{equation*}
\pi_{\mathfrak{m},\,\mathfrak{a}}:C(\mathfrak{m})\rightarrow C(\mathfrak{a})
\end{equation*}
be the natural surjective homomorphism (cf. \cite[Lemma 6.1 in Chapter V]{Janusz}).
Then we deduce that for a character $\xi$ of $C(\mathfrak{m})$
\begin{equation}\label{xicriterion}
\begin{aligned}
&~\mathfrak{f}_\xi=\mathfrak{m}\\
\Longleftrightarrow &~\xi|_{\mathrm{ker}(\pi_{\mathfrak{m},\,\mathfrak{m}\mathfrak{p}^{-1}})}\neq1~\textrm{for every}~\mathfrak{p}\mid\mathfrak{m}\\
\Longleftrightarrow &~\textrm{for each $\mathfrak{p}\mid\mathfrak{m}$, there is 
an element $C_\mathfrak{p}\in \mathrm{ker}(\pi_{\mathfrak{m},\,\mathfrak{m}\mathfrak{p}^{-1}})$ such that
$\xi(C_\mathfrak{p})\neq1$}. 
\end{aligned}
\end{equation}
For each $\mathfrak{p}\mid\mathfrak{m}$, take an element $u_\mathfrak{p}$ of $(\mathcal{O}_K/\mathfrak{m})^\times$
which reduces to the identity
of $(\mathcal{O}_K/\mathfrak{m}\mathfrak{p}^{-1})^\times$ 
as in Table \ref{table2}. The last column of Table 
\ref{table2} follows directly from \eqref{x1x2x3}. 
\begin{center}
\begin{tabular}{@{}cccc@{}}
\toprule
$\mathfrak{m}$ & $\mathfrak{p}\mid\mathfrak{m}$ & $u_\mathfrak{p}$
in the coordinates of Table \ref{table1}
& $\theta(u_\mathfrak{p})$ \\
\midrule\midrule
$\mathfrak{p}_{3,\,\mathrm{split}}^2$ & $\mathfrak{p}_{3,\,\mathrm{split}}$
& $4$ & $\zeta_3^2$\\ 
\midrule
$\mathfrak{p}_{3,\,\mathrm{split}}^2\mathfrak{p}'_{3,\,\mathrm{split}}$ & 
$\mathfrak{p}_{3,\,\mathrm{split}}$ & $(4,\,1)$  & $\zeta_3$ \\ 
& $\mathfrak{p}'_{3,\,\mathrm{split}}$ & $(1,\,-1)$ & $-1$ \\ 
\midrule
$\mathfrak{p}_{3,\,\mathrm{noninert}}\mathfrak{p}_{5,\,\mathrm{noninert}}$ & 
$\mathfrak{p}_{3,\,\mathrm{noninert}}$ & $(-1,\,1)$ & $-1$ \\
 & $\mathfrak{p}_{5,\,\mathrm{noninert}}$ & $(1,\,2)$ & $\mathrm{i}$ \\ 
\midrule
$\mathfrak{p}_{3,\,\mathrm{split}}\mathfrak{p}'_{3,\,\mathrm{split}}\mathfrak{p}_{5,\,\mathrm{noninert}}$ & $\mathfrak{p}_{3,\,\mathrm{split}}$ & $(-1,\,1,\,1)$  & $-1$ \\
& $\mathfrak{p}'_{3,\,\mathrm{split}}$ & $(1,\,-1,\,1)$ & $-1$\\
& $\mathfrak{p}_{5,\,\mathrm{noninert}}$ & $(1,\,1,\,2)$ & $-1$\\  
\bottomrule
\end{tabular}
\captionof{table}{An element $u_\mathfrak{p}\in(\mathcal{O}_K/\mathfrak{m})^\times$
reducing to the identity of $(\mathcal{O}_K/\mathfrak{m}\mathfrak{p}^{-1})^\times$}\label{table2}
\end{center}
For each $\mathfrak{p}\mid\mathfrak{m}$, choose a representative $a_\mathfrak{p}\in\mathcal{O}_K$
of $u_\mathfrak{p}\in(\mathcal{O}_K/\mathfrak{m})^\times$. Then 
we observe that
$a_\mathfrak{p}\mathcal{O}_K$ is relatively prime to $\mathfrak{m}$ and 
\begin{equation*}
a_\mathfrak{p}\equiv1\Mod{\mathfrak{m}\mathfrak{p}^{-1}}. 
\end{equation*} 
And, the class
\begin{equation*}
A_\mathfrak{p}:=q_\mathfrak{m}(u_\mathfrak{p})=[a_\mathfrak{p}\mathcal{O}_K]_\mathfrak{m}\in R_\mathfrak{m}(=\mathrm{ker}(C(\mathfrak{m})\rightarrow C(\mathcal{O}_K)))
\end{equation*}
satisfies
\begin{equation}\label{APtrivial}
\pi_{\mathfrak{m},\,\mathfrak{m}\mathfrak{p}^{-1}}(A_\mathfrak{p})=[a_\mathfrak{p}\mathcal{O}_K]_{\mathfrak{m}\mathfrak{p}^{-1}}=[\mathcal{O}_K]_{\mathfrak{m}\mathfrak{p}^{-1}}
\end{equation}
and
\begin{equation}\label{paet}
\eta(A_\mathfrak{p})=
\widetilde{\theta}(q_\mathfrak{m}(u_\mathfrak{p}))=\theta(u_\mathfrak{p})\neq 1.
\end{equation}
Hence we obtain by \eqref{xicriterion} that
\begin{equation*}
\mathfrak{f}_\eta=\mathfrak{m}. 
\end{equation*}
Here, if $\eta|_{\sigma_\mathfrak{m}^{-1}(\mathrm{Gal}(K_\mathfrak{m}/L))}\neq1$, then $\chi=\eta$ satisfies the required property \eqref{required}, and we are done. 
\par
Otherwise, 
if
\begin{equation}\label{etatrivial}
\eta|_{\sigma_\mathfrak{m}^{-1}(\mathrm{Gal}(K_\mathfrak{m}/L))}=1,
\end{equation}
choose an element 
\begin{equation*}
A\in\sigma_\mathfrak{m}^{-1}(\mathrm{Gal}(K_\mathfrak{m}/L))\setminus\{[\mathcal{O}_K]_\mathfrak{m}\}
\end{equation*}
by \eqref{nontrivialL}. Since $A\not\in R_\mathfrak{m}$ by \eqref{intersection}, $AR_\mathfrak{m}$ in $C(\mathfrak{m})/R_\mathfrak{m}$ is nontrivial. So there is a character $\psi$ of the abelian group $C(\mathfrak{m})/R_\mathfrak{m}$
such that
\begin{equation}\label{psinontrivial}
\psi(AR_\mathfrak{m})\neq1
\end{equation}
(see \cite[p. 61]{Serre}). Let $\widetilde{\psi}$ be
its inflation to $C(\mathfrak{m})$ so that
\begin{equation}\label{psitilde}
\widetilde{\psi}|_{R_\mathfrak{m}}=1
\end{equation}
and
\begin{equation*}
\widetilde{\psi}(B)=\psi(BR_\mathfrak{m})\quad(B\in  C(\mathfrak{m})).
\end{equation*} 
Since $C(\mathfrak{m})/R_\mathfrak{m}\simeq C(\mathcal{O}_K)$, 
we attain $\mathfrak{f}_{\widetilde{\psi}}=\mathcal{O}_K$. Put
\begin{equation*}
\chi=\eta\widetilde{\psi}. 
\end{equation*}
We see from \eqref{etatrivial} and \eqref{psinontrivial} that
\begin{equation*}
\chi(A)=\eta(A)\widetilde{\psi}(A)=\psi(AR_\mathfrak{m})\neq1,
\end{equation*}
and so
\begin{equation*}
\chi|_{\sigma_\mathfrak{m}^{-1}(\mathrm{Gal}(K_\mathfrak{m}/L))}\neq1.
\end{equation*}
For each $\mathfrak{p}\mid\mathfrak{m}$, recall that $A_\mathfrak{p}\in R_\mathfrak{m}$, and observe by \eqref{psitilde} and
\eqref{paet} that
\begin{equation*}
\chi(A_\mathfrak{p})=\eta(A_\mathfrak{p})\widetilde{\psi}(A_\mathfrak{p})
=\eta(A_\mathfrak{p})\neq1.
\end{equation*}
Therefore \eqref{xicriterion} and \eqref{APtrivial} give $\mathfrak{f}_{\chi}=\mathfrak{m}$,
and hence $\chi$ is a required character of $C(\mathfrak{m})$ satisfying \eqref{required}.  

\end{proof}

\begin{theorem}\label{main1}
Let $K$ be an imaginary quadratic field,
and let $E$ and $h_E$ be as in \eqref{model}--\eqref{Weberh}. 
If $\mathfrak{m}$ is a nonzero proper ideal of $\mathcal{O}_K$
that is the exact conductor of $K_\mathfrak{m}/K$ and is one of four cases
\begin{equation}\label{4forms2}
\mathfrak{p}_{3,\,\mathrm{split}}^2,\quad
\mathfrak{p}_{3,\,\mathrm{split}}^2\mathfrak{p}'_{3,\,\mathrm{split}},\quad
\mathfrak{p}_{3,\,\mathrm{noninert}}\mathfrak{p}_{5,\,\mathrm{noninert}},\quad
\mathfrak{p}_{3,\,\mathrm{split}}\mathfrak{p}'_{3,\,\mathrm{split}}\mathfrak{p}_{5,\,\mathrm{noninert}},
\end{equation}
then we have
\begin{equation*}
K_\mathfrak{m}=K(h_E(P))
\end{equation*}
for every primitive $\mathfrak{m}$-division
point $P\in E(\mathbb{C})$.
\end{theorem}
\begin{proof}
If $h_K=1$, then the result follows directly from \eqref{overH}. 
\par
Now, consider the case where $h_K>1$. Let 
\begin{equation*}
L=K(h_E(P)).  
\end{equation*}
Since $K_\mathfrak{m}/K$ is an abelian extension, 
so is $L/K$. 
Suppose on the contrary that $L$ is properly contained in $K_\mathfrak{m}$. 
Then one can choose a character $\chi$ satisfying
\begin{equation}\label{required1}
\mathfrak{f}_\chi=\mathfrak{m}
\end{equation}
and
\begin{equation}\label{required2}
\chi|_{\sigma_\mathfrak{m}^{-1}(\mathrm{Gal}(K_\mathfrak{m}/L))}\neq1 
\end{equation}
by Lemma \ref{charlemma}. 
Since $2+\mathfrak{m}\in(\mathcal{O}_K/\mathfrak{m})^\times$ with 
$\mathfrak{m}$ listed in \eqref{4forms2}, we get by Lemma \ref{hconjugate} that
\begin{equation*}
h_E(2P)=h_E(P)^{\sigma_\mathfrak{m}([2\mathcal{O}_K]_\mathfrak{m})}\in L,
\end{equation*} 
so
\begin{equation*}
\delta:=h_E(2P)-h_E(P)\in L.
\end{equation*}
By Lemma \ref{Schilemma}, we have $\delta\neq0$. 
We then find that
\begin{align*}
&\sum_{C\in C(\mathfrak{m})}\chi(C)\log|(\delta^{12N_\mathfrak{m}})^{\sigma_\mathfrak{m}(C)}|\\
=&\sum_{
A\in\sigma_\mathfrak{m}^{-1}(\mathrm{Gal}(K_\mathfrak{m}/L))}
\sum_{B\in \mathcal{R}}
\chi(AB)\log|(\delta^{12N_\mathfrak{m}})^{\sigma_\mathfrak{m}(AB)}|\\
&\hspace{3cm}
\textrm{where $\mathcal{R}$ is a complete set of representatives
of}~
C(\mathfrak{m})/\sigma_\mathfrak{m}^{-1}(\mathrm{Gal}(K_\mathfrak{m}/L))\\
&=\left(\sum_{A\in\sigma_\mathfrak{m}^{-1}(\mathrm{Gal}(K_\mathfrak{m}/L))}\chi(A)\right)\cdot
\left(\sum_{B\in \mathcal{R}}
\chi(B)\log|(\delta^{12N_\mathfrak{m}})^{\sigma_\mathfrak{m}(B)}|\right)
\quad
\textrm{since}~
(\delta^{12N_\mathfrak{m}})^{\sigma_\mathfrak{m}(A)}=\delta^{12N_\mathfrak{m}}\\
&=0\cdot\left(\sum_{B\in \mathcal{R}}
\chi(B)\log|(\delta^{12N_\mathfrak{m}})^{\sigma_\mathfrak{m}(B)}|\right)\quad
\textrm{by \eqref{required2} and \cite[Proposition 4 in Chapter VI]{Serre}}\\
&=0.
\end{align*}
On the other hand, we deduce from 
\eqref{Lnonzero}, \eqref{required1}, Proposition \ref{LS} and Remark \ref{barremark} that
\begin{equation}\label{Snonzero}
S_\mathfrak{m}(\chi)\neq0. 
\end{equation}
We find that
\begin{align*}
\sum_{C\in C(\mathfrak{m})}\chi(C)\log|(\delta^{12N_\mathfrak{m}})^{\sigma_\mathfrak{m}(C)}|&
=-\overline{\chi}(C_P)(1+2\overline{\chi}([2\mathcal{O}_K]_\mathfrak{m}))S_\mathfrak{m}(\chi)
\quad\textrm{by Lemma \ref{Schilemma}}\\
&\neq0\quad\textrm{by \eqref{Snonzero}, since}~|\overline{\chi}([2\mathcal{O}_K]_\mathfrak{m})|=1,
\end{align*}
which yields a contradiction. 
Therefore we conclude that $L=K_\mathfrak{m}$, as desired. 
\end{proof}

\section {Proof by inequality arguments}

We treat the remaining three families by comparing the absolute
value of a difference of Weber values with 
those of its Galois conjugates. 

\begin{lemma}\label{jinequality}
If $h_K>1$ and $|d_K|\geq20$, then we have
\begin{equation*}
\left|\frac{j(\omega_Q)^2(j(\omega_Q)-1728)^3}{j(\tau_K)^2(j(\tau_K)-1728)^3}
\right|<877383\,|q_{\tau_K}|^{5/2}
\end{equation*}
for all nonprincipal reduced forms $Q\in\mathcal{Q}_\mathrm{red}(d_K)$. 
\end{lemma}
\begin{proof}
See \cite[Lemma 6.3 (ii)]{J-K-S}.
\end{proof}

\begin{lemma}\label{upperlower}
Let $(N,\,t)\in\{(12,\,5),~(8,\,3),~(5,\,2)\}$ and 
$\mathbf{w}=\begin{bmatrix}w_1&w_2\end{bmatrix}\in\frac{1}{N}M_{1,\,2}(\mathbb{Z})
\setminus M_{1,\,2}(\mathbb{Z})$.
\begin{enumerate}
\item[\textup{(i)}] If $\tau\in\mathbb{H}$ satisfies $\mathrm{Im}(\tau)\geq\sqrt{3}/{2}$, then
\begin{equation*}
\left|\wp(t(w_1\tau+w_2);\,[\tau,\,1])-\wp(w_1\tau+w_2;\,[\tau,\,1])\right|<\frac{(2\pi)^2N^2}{10}.
\end{equation*}
\item[\textup{(ii)}] If $|d_K|\geq15$, then for every integer $s$
\begin{equation*}
\left|
\wp\left(\frac{t(\tau_K+s)}{N};\,[\tau_K,\,1]\right)-\wp\left(\frac{\tau_K+s}{N};\,[\tau_K,\,1]\right)
\right|>\frac{(2\pi)^2}{2}|q_{\tau_K}|^{1/N}.
\end{equation*}
\end{enumerate}
\end{lemma}
\begin{proof}
By Lemma \ref{pexpansion}, we have the expansion formula
\begin{equation}\label{precall}
\begin{aligned}
\frac{1}{(2\pi\mathrm{i})^2}\,\wp(z;\,[\tau,\,1])&=
\frac{1}{12}+\frac{q_z}{(1-q_z)^2}\\
&+\sum_{n=1}^\infty
\left(
\frac{q_\tau^nq_z}{(1-q_\tau^nq_z)^2}+
\frac{q_\tau^nq_z^{-1}}{(1-q_\tau^nq_z^{-1})^2}\right)
-2\sum_{n=1}^\infty\frac{nq_\tau^n}{1-q_\tau^n}
\quad(\tau\in\mathbb{H},~z\in\mathbb{C}).
\end{aligned}
\end{equation}
\begin{enumerate}
\item[(i)] Let $\tau\in\mathbb{H}$ with $\mathrm{Im}(\tau)\geq\sqrt{3}/2$. 
Take $\mathbf{u}=\begin{bmatrix}u_1&u_2\end{bmatrix},\,
\mathbf{v}=\begin{bmatrix}v_1&v_2\end{bmatrix}\in\frac{1}{N}M_{1,\,2}(\mathbb{Z})
\setminus M_{1,\,2}(\mathbb{Z})$ such that
\begin{equation}\label{uv}
\left\{\begin{array}{l}
0\leq u_1\leq\displaystyle\frac{1}{2}
\quad\textrm{and}\quad
\mathbf{u}\equiv\mathbf{w}~
\textrm{or}~
-\mathbf{w}\Mod{M_{1,\,2}(\mathbb{Z})},\vspace{0.1cm}\\
0\leq v_1\leq\displaystyle\frac{1}{2}
\quad\textrm{and}\quad
\mathbf{v}\equiv t\mathbf{w}~
\textrm{or}~
-t\mathbf{w}\Mod{M_{1,\,2}(\mathbb{Z})}.
\end{array}\right.
\end{equation}
Put $\rho:=|q_\tau|$. Since $\tau\in\mathbb{H}$, we see that $0<\rho<1$.
We then obtain by \eqref{precall} and \eqref{uv} that
\begin{equation}\label{differencew}
\begin{aligned}
&\frac{1}{(2\pi)^2}\left|\wp(t(w_1\tau+w_2);\,[\tau,\,1])-\wp(w_1\tau+w_2;\,[\tau,\,1])\right|\\
=&\frac{1}{(2\pi)^2}\left|\wp(v_1\tau+v_2;\,[\tau,\,1])-\wp(u_1\tau+u_2;\,[\tau,\,1])\right|\\
&\textrm{by evenness and periodicity of the Weierstrass $\wp$-function (see \cite[Lemma 10.4]{Cox})}\\\leq&
\left|\frac{q_{z'}}{(1-q_{z'})^2}\right|
+\left|\frac{q_z}{(1-q_z)^2}\right|\\
&+
\sum_{n=1}^\infty
\left(
\frac{\rho^{n+v_1}}{(1-\rho^{n+v_1})^2}
+\frac{\rho^{n-v_1}}{(1-\rho^{n-v_1})^2}\right)
+
\sum_{n=1}^\infty
\left(
\frac{\rho^{n+u_1}}{(1-\rho^{n+u_1})^2}
+\frac{\rho^{n-u_1}}{(1-\rho^{n-u_1})^2}\right)\\
&\hspace{7cm}\textrm{where}~z'=v_1\tau+v_2~
\textrm{and}~z=u_1\tau+u_2\\
\leq&
\left|\frac{q_{z'}}{(1-q_{z'})^2}\right|
+\left|\frac{q_z}{(1-q_z)^2}\right|
+2
\sum_{n=1}^\infty
\left(
\frac{\rho^n}{(1-\rho^n)^2}
+\frac{\rho^{n-1/2}}{(1-\rho^{n-1/2})^2}\right).
\end{aligned}
\end{equation}
Observe that if $u_1=0$, then 
\begin{align*}
\left|\frac{q_z}{(1-q_z)^2}\right|&=
\frac{1}{4\sin^2(\pi u_2)}\\
&\leq\frac{1}{4\sin^2(\pi/N)}\\
&\leq\frac{N^2}{32}\quad\textrm{due to}~
\sin\,x\geq\frac{2\sqrt{2}}{\pi}x~(0\leq x\leq\frac{\pi}{4})~
\textrm{and}~N>4,
\end{align*}
and if $u_1\neq0$, then 
\begin{align*}
\left|\frac{q_z}{(1-q_z)^2}\right|&\leq
\frac{|q_z|}{(1-|q_z|)^2}\\
&=\frac{1}{4\,\mathrm{sinh}^2(\pi y)}\quad\textrm{where}~y=\mathrm{Im}(z)\\
&<\frac{1}{(2\pi y)^2}\quad\textrm{since}~\mathrm{sinh}\,x>x~(x>0)\\
&\leq\frac{N^2}{3\pi^2}\quad\textrm{because}~y=u_1\mathrm{Im}(\tau)\geq\frac{\sqrt{3}}{2N}\\
&<\frac{N^2}{25}.
\end{align*}
Hence it follows that
\begin{equation}\label{firsttwo}
\left|\frac{q_z}{(1-q_z)^2}\right|<\frac{N^2}{25},\quad\textrm{similarly}\quad
\left|\frac{q_{z'}}{(1-q_{z'})^2}\right|<\frac{N^2}{25}.
\end{equation}
Furthermore, we find that
\begin{equation}\label{sum2}
\begin{aligned}
\sum_{n=1}^\infty\left(
\frac{\rho^n}{(1-\rho^n)^2}+
\frac{\rho^{n-1/2}}{(1-\rho^{n-1/2})^2}\right)
&\leq\sum_{n=1}^\infty\left(
\frac{\rho^n}{(1-\rho)^2}+\frac{\rho^{n-1/2}}{(1-\rho^{1/2})^2}
\right)\\
&\leq\frac{\rho}{(1-\rho)^3}+\frac{\rho^{1/2}}
{(1-\rho)(1-\rho^{1/2})^2}\\
&<\frac{1}{10}\quad\textrm{because}~\mathrm{Im}(\tau)\geq\sqrt{3}/{2}\\
&<\frac{N^2}{100}\quad\textrm{since}~N\geq5. 
\end{aligned}
\end{equation}
Therefore, it follows from \eqref{differencew}--\eqref{sum2} that
\begin{equation*}
\left|\wp(t(w_1\tau+w_2);\,[\tau,\,1])-\wp(w_1\tau+w_2;\,[\tau,\,1])\right|<\frac{(2\pi)^2N^2}{10}.
\end{equation*}
\item[(ii)] 
Put
\begin{equation*}
r:=|q_{\tau_K}|^{1/N}=e^{-\pi\sqrt{|d_K|}/N}\quad\textrm{and}\quad
z:=\frac{\tau_K+s}{N}. 
\end{equation*}
The assumption $|d_K|\geq15$ implies
\begin{equation}\label{restimate}
r\leq e^{-\pi\sqrt{15}/N}<\left\{\begin{array}{ll}
3/8 & \textrm{if}~N=12,\\
1/4 & \textrm{if}~N=8,\\
1/10 & \textrm{if}~N=5. 
\end{array}\right.
\end{equation}
Since $0<1/N<t/N<1$, we see that neither $z$ nor $tz$ belongs to $[\tau_K,\,1]$. 
We then derive by \eqref{precall} that
\begin{equation}\label{2pipdifference}
\begin{aligned}
\frac{1}{(2\pi)^2}\left|\wp(tz;\,[\tau_K,\,1])-
\wp(z;\,[\tau_K,\,1])\right|
&\geq 
\left|\frac{q_z}{(1-q_z)^2}\right|-
\left|\frac{q_{tz}}{(1-q_{tz})^2}\right|-T\\
&\geq\frac{r}{(1+r)^2}-\frac{r^t}{(1-r^t)^2}-T,
\end{aligned}
\end{equation}
where
\begin{equation*}
T=
\sum_{n=1}^\infty\left(
\frac{r^{Nn+1}}{(1-r^{Nn+1})^2}+\frac{r^{Nn-1}}{(1-r^{Nn-1})^2}
+\frac{r^{Nn+t}}{(1-r^{Nn+t})^2}+\frac{r^{Nn-t}}{(1-r^{Nn-t})^2}
\right).
\end{equation*}
Meanwhile we also find that
\begin{equation}\label{Tinequality}
\begin{aligned}
T&\leq\frac{1}{(1-r^{N-t})^2}\left\{\sum_{n=1}^\infty r^{Nn+1}
+\sum_{n=1}^\infty r^{Nn-1}
+\sum_{n=1}^\infty r^{Nn+t}
+\sum_{n=1}^\infty r^{Nn-t}\right\}\\
&=\frac{r^{N+1}+r^{N-1}+r^{N+t}+r^{N-t}}{(1-r^{N-t})^2(1-r^N)}\\
&=\frac{r^{N-t}(1+r^{t-1})(1+r^{t+1})}{(1-r^{N-t})^2(1-r^N)}\\
&<\frac{2r^{N-t}}{{(1-r^{N-t})^2(1-r^N)}}\quad\textrm{
since}~(1+r^{t-1})(1+r^{t+1})<2\\
&\hspace{4.5cm}\textrm{owing to \eqref{restimate} and the fact $t\geq2$}. 
\end{aligned}
\end{equation}
Thus we achieve by \eqref{2pipdifference} and \eqref{Tinequality} that
\begin{align*}
&\frac{1}{(2\pi)^2}\left|\wp(tz;\,[\tau_K,\,1])-
\wp(z;\,[\tau_K,\,1])\right|\\
>&~r\left(\frac{1}{(1+r)^2}-\frac{r^{t-1}}{(1-r^t)^2}-
\frac{2r^{N-t-1}}{{(1-r^{N-t})^2(1-r^N)}}\right)\\
>&~r\cdot\frac{1}{2}\quad\textrm{by \eqref{restimate} and the fact}~
(N,\,t)\in\{(12,\,5),~(8,\,3),~(5,\,2)\},
\end{align*}
and so
\begin{equation*}
\left|
\wp\left(\frac{t(\tau_K+s)}{N};\,[\tau_K,\,1]\right)-\wp\left(\frac{\tau_K+s}{N};\,[\tau_K,\,1]\right)
\right|>\frac{(2\pi)^2}{2}|q_{\tau_K}|^{1/N}.
\end{equation*}
\end{enumerate}
\end{proof}

\begin{theorem}\label{main2}
Let $K$ be an imaginary quadratic field,
and let $E$ and $h_E$ be as in \eqref{model}--\eqref{Weberh}. 
If $\mathfrak{m}$ is a nonzero proper ideal of $\mathcal{O}_K$
that is the exact conductor of $K_\mathfrak{m}/K$ and is one of 
three cases
\begin{equation*}
\mathfrak{p}_{2,\,\mathrm{split}}^2\mathfrak{p}_{3,\,\mathrm{noninert}},\quad
\mathfrak{p}_{2,\,\mathrm{split}}^3,\quad
\mathfrak{p}_{5,\,\mathrm{noninert}},
\end{equation*}
then we have
\begin{equation*}
K_\mathfrak{m}=K(h_E(P))
\end{equation*}
for every primitive $\mathfrak{m}$-division
point $P\in E(\mathbb{C})$.
\end{theorem}
\begin{proof}
If $h_K=1$, then the assertion follows from \eqref{overH}. So we assume that 
\begin{equation}\label{15}
h_K>1,~\textrm{and so}~|d_K|\geq15~\textrm{and}~\mathcal{O}_K^\times=\{\pm1\}
\end{equation}
by \eqref{classnumberone} and \eqref{O_Ktimes}. Set
\begin{equation*}
(N,\,t)=\left\{\begin{array}{ll}
(12,\,5) & \textrm{if}~\mathfrak{m}=\mathfrak{p}_{2,\,\mathrm{split}}^2\mathfrak{p}_{3,\,\mathrm{noninert}},\\
(8,\,3) & \textrm{if}~\mathfrak{m}=\mathfrak{p}_{2,\,\mathrm{split}}^3,\\
(5,\,2) & \textrm{if}~\mathfrak{m}=\mathfrak{p}_{5,\,\mathrm{noninert}}.
\end{array}\right.
\end{equation*}
Since $\mathrm{N}(\mathfrak{m})=N$ and 
$\mathfrak{m}\cap\mathbb{Z}=N\mathbb{Z}$,
it follows from \cite[Theorem 6.9 and Proposition 6.11]{Buell} that
\begin{equation*}
\overline{\mathfrak m}=[\tau_K+s,\,N]\quad\textrm{for some}~s\in\mathbb{Z}, 
\end{equation*}
and so
\begin{equation*}
\mathfrak{m}^{-1}=\frac{1}{N}\overline{\mathfrak{m}}=[z,\,1]\quad
\textrm{with}~z=\frac{\tau_K+s}{N}. 
\end{equation*}
We then find that
\begin{equation*}
\mathrm{Ann}_{\mathcal{O}_K}(\phi_E(z+\mathcal{O}_K))=\{a\in\mathcal{O}_K\mid
az\in\mathcal{O}_K\}=\{a\in\mathcal{O}_K\mid a\mathfrak{m}^{-1}\subseteq\mathcal{O}_K\}
=\mathfrak{m}, 
\end{equation*}
which claims that $\phi_E(z+\mathcal{O}_K)$ is a primitive $\mathfrak{m}$-division point. 
Since the Weber values at primitive $\mathfrak{m}$-division points are Galois conjugates
over $H_K$, it suffices to prove the assertion for $P=\phi_E(z+\mathcal{O}_K)$.
\par
Let $P=\phi_E(z+\mathcal{O}_K)$ and $L=K(h_E(P))$. Since $K_\mathfrak{m}/K$
is an abelian extension, so is $L/K$. 
Suppose on the contrary that $L$ is properly contained in $K_\mathfrak{m}$. 
Since $(t,\,N)=1$, we get by Lemma \ref{hconjugate} that
\begin{equation*}
h_E(tP)=h_E(P)^{\sigma_\mathfrak{m}([t\mathcal{O}_K]_\mathfrak{m})}\in L,
\end{equation*}
and so
\begin{equation}\label{differenceL}
\delta:=h_E(tP)-h_E(P)\in L. 
\end{equation}
Here, we note by \eqref{jdef}, \eqref{not01728}, \eqref{15} and Lemma \ref{upperlower} (ii) that $\delta\neq0$. 
Take an element 
\begin{equation*}
\varepsilon\in\mathrm{Gal}(K_\mathfrak{m}/L)\setminus\{\mathrm{id}_{K_\mathfrak{m}}\}, 
\end{equation*}
and hence
\begin{equation}\label{delta=deltae}
\delta^\varepsilon=\delta
\end{equation}
by \eqref{differenceL}. Furthermore, since
$H_KL=K_\mathfrak{m}$ by \eqref{overH}, we see that
\begin{equation}\label{notid}
\varepsilon|_{H_K}\neq
\mathrm{id}|_{H_K}. 
\end{equation}
Let $C=\sigma_\mathfrak{m}^{-1}(\varepsilon)$, and 
choose an ideal $\mathfrak{c}$ of $\mathcal{O}_K$ in $C$. 
Then we find that 
\begin{equation}\label{he}
\begin{aligned}
h_E(P)^\varepsilon&=f_\mathfrak{m}(C_P)^{\sigma_\mathfrak{m}(C)}\quad\textrm{by Remark \ref{CP} (i)}\\
&=f_\mathfrak{m}(C_PC)\quad\textrm{by Proposition \ref{fginvariants}},
\end{aligned}
\end{equation}
and similarly
\begin{equation}\label{hte}
\begin{aligned}
h_E(tP)^\varepsilon&=\left(h_E(P)^{\sigma_\mathfrak{m}([t\mathcal{O}_K]_\mathfrak{m})}\right)
^{\sigma_\mathfrak{m}(C)}
\quad\textrm{by Lemma \ref{hconjugate}}\\
&=f_\mathfrak{m}(C_P)^{\sigma_\mathfrak{m}([t\mathcal{O}_K]_\mathfrak{m}C)}
\quad\textrm{by Remark \ref{CP} (i)}\\
&=f_\mathfrak{m}([t\mathcal{O}_K]_\mathfrak{m}C_PC)\quad\textrm{by Proposition \ref{fginvariants}}.
\end{aligned}
\end{equation}
Observe that $C_PC=[z\mathfrak{m}\mathfrak{c}]_\mathfrak{m}$, and
\begin{equation*}
\mathfrak{m}(z\mathfrak{m}\mathfrak{c})^{-1}=z^{-1}\mathfrak{c}^{-1}=
\nu[\omega_Q,\,1]
\end{equation*}
for some $\nu\in K^\times$ and
$Q=ax^2+bxy+cy^2\in\mathcal{Q}_\mathrm{red}(d_K)$ by \eqref{bijective}. 
Here, $Q$ is nonprincipal by \eqref{notid}. 
Since $z\mathfrak{m}\mathfrak{c}$ is an ideal of $\mathcal{O}_K$ and $N\in\mathfrak{m}$, 
there is a vector $\mathbf{w}=\begin{bmatrix}w_1&w_2\end{bmatrix}\in
\frac{1}{N}M_{1,\,2}(\mathbb{Z})\setminus M_{1,\,2}(\mathbb{Z})$ such that
\begin{equation*}
1=w_1(\nu\omega_Q)+w_2\nu. 
\end{equation*}
Thus Definition \ref{invariants} yields
\begin{equation*}
f_\mathfrak{m}(C_PC)=f_\mathbf{w}(\omega_Q). 
\end{equation*}
In a similar way, one derives
\begin{equation*}
f_\mathfrak{m}([t\mathcal{O}_K]_\mathfrak{m}C_PC)=f_{t\mathbf{w}}(\omega_Q),
\end{equation*}
and hence
\begin{equation}\label{deltae}
\delta^\varepsilon=f_{t\mathbf{w}}(\omega_Q)-f_\mathbf{w}(\omega_Q) 
\end{equation}
by \eqref{he} and \eqref{hte}. 
\par
Put
\begin{equation*}
z_Q:=w_1\omega_Q+w_2,\quad
\rho=|q_{\tau_K}|=e^{-\pi\sqrt{|d_K|}},\quad
\rho_Q:=|q_{\omega_Q}|=\rho^{1/a}.
\end{equation*}
Since  $Q\in\mathcal{Q}_\mathrm{red}(d_K)$ is a nonprincipal reduced form, 
we attain by \eqref{abound} that
\begin{equation}\label{2200}
a\geq2\quad\textrm{and}
\quad0<\rho<\rho_Q\leq e^{-\pi\sqrt{3}}<\frac{1}{200}. 
\end{equation}
We find that
\begin{equation}\label{|delta|inequality}
\begin{aligned}
|\delta|&=\left|
\frac{g_2(\tau_K)g_3(\tau_K)}{\Delta(\tau_K)}\right|
\left|\wp(tz;\,[\tau_K,\,1])-\wp(z;\,[\tau_K,\,1])\right|\\
&\hspace{5cm}\textrm{by the definitions \eqref{differenceL}, 
\eqref{Weberh} and \eqref{Frickedef}}\\
&>\frac{(2\pi)^2}{2}\left|\frac{g_2(\tau_K)g_3(\tau_K)}
{\Delta(\tau_K)}\right|\rho^{1/N}\quad\textrm{by Lemma \ref{upperlower} (ii)}\\
&=\frac{(2\pi)^2}{2}\left|
\frac{j(\tau_K)^2(j(\tau_K)-1728)^3}
{2^{30}3^{24}\Delta(\tau_K)}
\right|^{1/6}\rho^{1/N}\quad\textrm{by the definitions \eqref{Deltadef} and \eqref{jdef}}
\end{aligned}
\end{equation}
and
\begin{equation}\label{|deltae|inequality}
\begin{aligned}
|\delta^\varepsilon|&=\left|
\frac{g_2(\omega_Q)g_3(\omega_Q)}{\Delta(\omega_Q)}\right|
\left|\wp(tz_Q;\,[\omega_Q,\,1])-\wp(z_Q;\,[\omega_Q,\,1])\right|\\
&\hspace{5cm}\textrm{by \eqref{deltae}, the definitions \eqref{Weberh} and \eqref{Frickedef}}\\
&<\frac{(2\pi)^2N^2}{10}\left|\frac{g_2(\omega_Q)g_3(\omega_Q)}
{\Delta(\omega_Q)}\right|\quad\textrm{by Lemma \ref{upperlower} (i)}\\
&=\frac{(2\pi)^2N^2}{10}
\left|
\frac{j(\omega_Q)^2(j(\omega_Q)-1728)^3}
{2^{30}3^{24}\Delta(\omega_Q)}
\right|^{1/6}\quad\textrm{by the definitions \eqref{Deltadef} and \eqref{jdef}}.
\end{aligned}
\end{equation}
Furthermore, we obtain that
\begin{equation}\label{Deltaquotient}
\begin{aligned}
\left|
\frac{\Delta(\tau_K)}{\Delta(\omega_Q)}\right|
&\leq\frac{\rho}{\rho_Q}
\prod_{n=1}^\infty
\left(\frac{1+\rho^n}{1-\rho_Q^n}
\right)^{24}\quad\textrm{by \eqref{Deltaexpansion} of Lemma \ref{g2g3Deltaexpansion}}\\
&\leq\frac{\rho}{\rho_Q}
\prod_{n=1}^\infty
\left(\frac{1}{(1-\rho^n)(1-\rho_Q^n)}
\right)^{24}\quad\textrm{since}~1+x\leq\frac{1}{1-x}~(0\leq x<1)\\
&\leq\frac{\rho}{\rho_Q}
\left(\frac{1-\rho}{1-2\rho}\frac{1-\rho_Q}{1-2\rho_Q}\right)^{24}\\
&\hspace{2cm}\textrm{because}~\prod_{n=1}^\infty(1-x_n)\geq 1-\sum_{n=1}^\infty x_n~(x_n\geq0~
\textrm{with}~\sum_{n=1}^\infty x_n<1)\\
&<\rho^{1/2}\left(\frac{199}{198}\right)^{48}\quad\textrm{by \eqref{2200}}\\
&<2\rho^{1/2}.
\end{aligned}
\end{equation}
Now, we achieve by \eqref{delta=deltae}, \eqref{|delta|inequality},
\eqref{|deltae|inequality} and \eqref{Deltaquotient} that
\begin{equation}\label{inequalitycriterion}
1=\left|
\frac{\delta^\varepsilon}{\delta}
\right|^6<2\left(\frac{N^2}{5}\right)^6 \rho^{(1/2)-(6/N)}
\left|\frac{j(\omega_Q)^2(j(\omega_Q)-1728)^3}
{j(\tau_K)^2(j(\tau_K)-1728)^3}\right|.
\end{equation}

\begin{enumerate}
\item[Case 1.] Suppose that $|d_K|\geq20$. We then find that
\begin{align*}
1&<2\left(\frac{N^2}{5}\right)^6 \rho^{(1/2)-(6/N)}
\left|\frac{j(\omega_Q)^2(j(\omega_Q)-1728)^3}
{j(\tau_K)^2(j(\tau_K)-1728)^3}\right|
\quad\textrm{by \eqref{inequalitycriterion}}\\
&<2\left(\frac{N^2}{5}\right)^6 \rho^{(1/2)-(6/N)}\cdot
877383\,\rho^{5/2}\quad\textrm{by Lemma \ref{jinequality}}\\
&=1754766\left(\frac{N^2}{5}\right)^6\rho^{3-(6/N)}\\
&\leq1754766\left(\frac{N^2}{5}\right)^6(e^{-\pi\sqrt{20}})^{3-(6/N)}\quad\textrm{due to}~|d_K|\geq20\\
&<1\quad\textrm{for every}~N=12,\,8,\,5,
\end{align*}
which yields a contradiction. 
\item[Case 2.] Suppose that $d_K=-15$. In this case, we have
$\mathcal{Q}_\mathrm{red}(d_K)=\{x^2+xy+4y^2,\,2x^2+xy+2y^2\}$,
in which $Q=2x^2+xy+2y^2$ is the unique nonprincipal form.
Moreover, we have
\begin{equation*}
j(\tau_K)=\frac{-191025-85995\sqrt{5}}{2}\quad\textrm{and}\quad
j(\omega_Q)=\frac{-191025+85995\sqrt{5}}{2} 
\end{equation*}
(see \cite[Example 6.2.2 in Chapter II]{Silverman94}), and estimate that
\begin{equation}\label{sharp}
\left|\frac{j(\omega_Q)^2(j(\omega_Q)-1728)^3}
{j(\tau_K)^2(j(\tau_K)-1728)^3}\right|<2\cdot10^{-12}.
\end{equation}
So, we establish that
\begin{align*}
1&<2\left(\frac{N^2}{5}\right)^6 \rho^{(1/2)-(6/N)}
\left|\frac{j(\omega_Q)^2(j(\omega_Q)-1728)^3}
{j(\tau_K)^2(j(\tau_K)-1728)^3}\right|\quad\textrm{by \eqref{inequalitycriterion}}\\
&<2\left(\frac{N^2}{5}\right)^6  (e^{-\pi\sqrt{15}})^{(1/2)-(6/N)}(2\cdot10^{-12})
\quad\textrm{by \eqref{sharp}}\\
&<1\quad\textrm{for every}~N=12,\,8,\,5,
\end{align*}
which yields a contradiction too. 
\end{enumerate}
Therefore we conclude that $L=K_\mathfrak{m}$ as desired. 
\end{proof}

Combining the preceding results with Proposition \ref{residual} yields 
Theorem \ref{main}. 

\begin{theorem}\label{conclusion}
Let $K$ be an imaginary quadratic field,
and let $E$ and $h_E$ be as in \eqref{model}--\eqref{Weberh}. 
If $\mathfrak{m}$ is a nonzero proper ideal of $\mathcal{O}_K$
that is the exact conductor of $K_\mathfrak{m}/K$, then
\begin{equation*}
K_\mathfrak{m}=K(h_E(P))
\end{equation*}
for every primitive $\mathfrak{m}$-division
point $P\in E(\mathbb{C})$.
\end{theorem}
\begin{proof}
The case $h_K=1$ follows from \eqref{overH}. For $h_K>1$, 
Sugawara's criteria (S1)--(S3) and \cite{J-K-S, K-S-Y} treat
the cases excluded in Proposition \ref{residual}. And, the remaining seven families
are covered by Theorems \ref{main1} and \ref{main2}. 
\end{proof}

\section {Generation of intermediate fields using relative norms}\label{sect8}

Starting from a Weber value that generates $K_\mathfrak{m}$ over $K$,
we construct an explicit integer translate whose relative norms
generate every intermediate field. 

For a nonzero proper ideal $\mathfrak{m}$ of $\mathcal{O}_K$, let $N_\mathfrak{m}>1$ be the integer
satisfying 
\begin{equation*}
\mathfrak{m}\cap\mathbb{Z}=N_\mathfrak{m}\mathbb{Z}.
\end{equation*}
We further define 
two integers $\mathcal{A}_\mathfrak{m}$ and $\mathcal{B}_\mathfrak{m}$ by 
\begin{equation*}
\mathcal{A}_\mathfrak{m}=2^93^6N_\mathfrak{m}^6\quad
\textrm{and}\quad
\mathcal{B}_\mathfrak{m}=\lceil N_\mathfrak{m}^6 e^{\pi\sqrt{|d_K|}}\rceil.
\end{equation*}

\begin{lemma}\label{integral}
Let $N>1$ be an integer and $\mathbf{v}\in\frac{1}{N}M_{1,\,2}(\mathbb{Z})
\setminus M_{1,\,2}(\mathbb{Z})$. Then the functions
\begin{equation*}
2^73^5N^2f^{(1)}_\mathbf{v},\quad
2^83^4N^4f^{(2)}_\mathbf{v},\quad
2^93^6N^6f^{(3)}_\mathbf{v}
\end{equation*}
are all integral over $\mathbb{Z}[j]$.
\end{lemma}
\begin{proof}

It is well known that
\begin{equation}\label{f1}
\begin{array}{l}
\text{the Fourier coefficients
of $2^73^5N^2f_\mathbf{v}^{(1)}\circ\gamma$ with respect to $q_\tau^{1/N}$}\\
\textrm{are algebraic integers
for every $\gamma\in\mathrm{SL}_2(\mathbb{Z})$}
\end{array}
\end{equation}
(see \cite[Lemma 5.6]{K-S} and Proposition \ref{modularGalois} (iii)). 
By the definitions \eqref{jdef} and \eqref{Frickedef}, we find that
\begin{equation}\label{f23}
2^83^4N^4f^{(2)}_\mathbf{v}=\frac{(2^73^5N^2f^{(1)}_\mathbf{v})^2}{j-1728}
\quad\textrm{and}\quad
2^93^6N^6f^{(3)}_\mathbf{v}=\frac{(2^73^5N^2f^{(1)}_\mathbf{v})^3}{j(j-1728)}.
\end{equation}
Observe from \eqref{jFourier} that
\begin{equation*}
q_\tau j,~q_\tau(j-1728)\in 1+q_\tau\mathbb{Z}[[q_\tau]],
\end{equation*}
and so
\begin{equation}\label{qq}
\frac{1}{j-1728}\in q_\tau\mathbb{Z}[[q_\tau]]\quad
\textrm{and}\quad
\frac{1}{j(j-1728)}\in q_\tau^2\mathbb{Z}[[q_\tau]]. 
\end{equation}
Thus we get by the fact
\eqref{f1}, the relations \eqref{f23} and 
\eqref{qq} that
\begin{equation}\label{f23gamma}
\begin{array}{l}
\text{the Fourier coefficients
of $2^83^4N^4f^{(2)}_\mathbf{v}\circ\gamma$ 
and $2^93^6N^6f^{(3)}_\mathbf{v}\circ\gamma$ with respect to $q_\tau^{1/N}$}\\
\textrm{are algebraic integers for every $\gamma\in\mathrm{SL}_2(\mathbb{Z})$}.
\end{array}
\end{equation}
Therefore, the fact from \cite[Lemma 2.1 in Chapter 2]{K-L}
together with \eqref{f1}, \eqref{f23gamma} and the holomorphy of 
$f^{(k)}_\mathbf{v}$ ($k=1,\,2,\,3$) on $\mathbb{H}$
lead us to conclude that
the functions
$2^73^5N^2f^{(1)}_\mathbf{v}$, $2^83^4N^4f^{(2)}_\mathbf{v}$ and $2^93^6N^6f^{(3)}_\mathbf{v}$ are  
integral over $\mathbb{Z}[j]$. 
\end{proof}

For every positive integer $N$, we get by the theory of complex multiplication that
\begin{equation*}
K_{N\mathcal{O}_K}=K(f(\tau_K)\mid f\in\mathcal{F}_N~\textrm{is finite at}~\tau_K)
\end{equation*}
(see \cite[Corollary to Theorem 2 in Chapter 10]{Lang}). 

\begin{lemma}\label{Aintegral}
If $P$ is a primitive $\mathfrak{m}$-division point of $E(\mathbb{C})$, then
$\mathcal{A}_\mathfrak{m}h_E(P)$ is an algebraic integer. 
\end{lemma}
\begin{proof}
Since $\mathcal{O}_K=[\tau_K,\,1]$ by \eqref{O_K} and $N_\mathfrak{m}\in\mathfrak{m}=\mathrm{Ann}_{\mathcal{O}_K}(P)$, 
we have
\begin{equation*}
\phi_E^{-1}(P)=(v_1\tau_K+v_2)+\mathcal{O}_K\quad
\textrm{for some}~\mathbf{v}=\begin{bmatrix}v_1&v_2\end{bmatrix}\in
\frac{1}{N_\mathfrak{m}}M_{1,\,2}(\mathbb{Z})\setminus M_{1,\,2}(\mathbb{Z}). 
\end{equation*}
And, we see from the definitions \eqref{Weberh} and
\eqref{Frickedef} that
\begin{equation*}
\mathcal{A}_\mathfrak{m}h_E(P)=\left\{
\begin{array}{ll}
2^23N_\mathfrak{m}^4\left(2^73^5N_\mathfrak{m}^2f^{(1)}_\mathbf{v}(\tau_K)\right) & \textrm{if}~j(E)\neq1728,\,0\vspace{0.1cm}\\
2^13^2N_\mathfrak{m}^2\left(2^83^4N_\mathfrak{m}^4f^{(2)}_\mathbf{v}(\tau_K)\right) & \textrm{if}~j(E)=1728,\vspace{0.1cm}\\
2^03^0N_\mathfrak{m}^0\left(2^93^6N_\mathfrak{m}^6f^{(3)}_\mathbf{v}(\tau_K)\right) & \textrm{if}~j(E)=0.
\end{array}\right.
\end{equation*}
Hence the value $\mathcal{A}_\mathfrak{m}h_E(P)$ is an algebraic integer
by Lemma \ref{integral} and the 
well-known fact that the singular modulus $j(\tau_K)$ is an algebraic integer (see \cite[Theorem 4 in Chapter 5]{Lang}). 
\end{proof}

\begin{lemma}\label{hinequality}
If $P$ is a primitive $\mathfrak{m}$-division point of $E(\mathbb{C})$, then
\begin{equation*}
|\iota(h_E(P))|\leq \mathcal{B}_\mathfrak{m}\quad
\textrm{for every embedding}~\iota:K_\mathfrak{m}\hookrightarrow\mathbb{C}. 
\end{equation*}
\end{lemma}
\begin{proof}
Put $k={|\mathcal{O}_K^\times|}/2$, and so
\begin{equation*}
k=\left\{\begin{array}{ll}
1 & \textrm{if}~K\neq\mathbb{Q}(\sqrt{-1}),\,\mathbb{Q}(\sqrt{-3}),\\
2 & \textrm{if}~K=\mathbb{Q}(\sqrt{-1}),\\
3 & \textrm{if}~K=\mathbb{Q}(\sqrt{-3})\\
\end{array}\right.
\end{equation*}
by \eqref{O_Ktimes}. 
Let $\iota:K_\mathfrak{m}\rightarrow\mathbb{C}$ be an embedding. 
Since $K$ is imaginary quadratic, $\iota|_K$ is either the identity or the complex conjugation. 
In the second case, the map 
\begin{equation*}
\iota':K_\mathfrak{m}\rightarrow\mathbb{C},\quad x\mapsto\overline{\iota(x)}
\end{equation*}
fixes $K$ and 
\begin{equation*}
|\iota(h_E(P))|=|\overline{\iota(h_E(P))}|=|\iota'(h_E(P))|.
\end{equation*}
Thus we may assume that $\iota$ fixes $K$. 
By Lemma \ref{fenough}, it suffices to show that
\begin{equation*}
|f_\mathbf{w}^{(k)}(\omega_Q)|\leq \mathcal{B}_\mathfrak{m}
\end{equation*}
for all $\mathbf{w}=\begin{bmatrix}w_1 & w_2\end{bmatrix}\in\frac{1}{N_\mathfrak{m}}M_{1,\,2}(\mathbb{Z})\setminus
M_{1,\,2}(\mathbb{Z})$ with $0\leq w_1\leq1/2$
and for all $Q\in\mathcal{Q}_\mathrm{red}(d_K)$.
\par 
By Lemma \ref{pexpansion}, we obtain the expansion formula
\begin{equation}\label{recallp}
\begin{aligned}
\frac{1}{(2\pi\mathrm{i})^2}\,\wp(z;\,[\tau,\,1])&=
\frac{1}{12}+\frac{q_z}{(1-q_z)^2}\\
&+\sum_{n=1}^\infty
\left(
\frac{q_\tau^nq_z}{(1-q_\tau^nq_z)^2}+
\frac{q_\tau^nq_z^{-1}}{(1-q_\tau^nq_z^{-1})^2}\right)
-2\sum_{n=1}^\infty\frac{nq_\tau^n}{1-q_\tau^n}
\quad(z\in\mathbb{C},~\tau\in\mathbb{H}).
\end{aligned}
\end{equation}
Let $\mathbf{w}=\begin{bmatrix}w_1 & w_2\end{bmatrix}\in\frac{1}{N_\mathfrak{m}}M_{1,\,2}(\mathbb{Z})
\setminus M_{1,\,2}(\mathbb{Z})$ with $0\leq w_1\leq1/2$, and let $Q=ax^2+bxy+cy^2
\in\mathcal{Q}_\mathrm{red}(d_K)$. Set
\begin{equation*}
z_Q:=w_1\omega_Q+w_2\quad\textrm{and}\quad 
y_Q:=\mathrm{Im}(\omega_Q)=\frac{\sqrt{|d_K|}}{2a}. 
\end{equation*}
Then we know by \eqref{abound} that
\begin{equation}\label{yinequality}
y_Q\geq\frac{\sqrt{3}}{2}.
\end{equation}
If $w_1=0$, then we attain that $w_2\not\in\mathbb{Z}$ and 
\begin{equation*}
|1-q_{z_Q}|=2|\sin(\pi w_2)|\geq2\sin\frac{\pi}{N_\mathfrak{m}}\geq\frac{4}{N_\mathfrak{m}}
\quad\textrm{since}~\sin(\pi x)\geq2x~(0\leq x\leq\frac{1}{2}), 
\end{equation*}
and so
\begin{equation}\label{S20}
\left|\frac{q_{z_Q}}{(1-q_{z_Q})^2}\right|\leq\frac{N_\mathfrak{m}^2}{16}. 
\end{equation}
If $0<w_1\leq1/2$, then
\begin{equation}\label{S2not0}
\begin{aligned}
\left|\frac{q_{z_Q}}{(1-q_{z_Q})^2}\right|&\leq\frac{|q_{z_Q}|}{(1-|q_{z_Q}|)^2}\\&=\frac{1}{4\mathrm{sinh}^2(\pi y_Qw_1)}\\
&<\frac{1}{4(\pi y_Qw_1)^2}\quad\textrm{because}~\mathrm{sinh}\,x> x~(x>0)\\
&<\frac{N_\mathfrak{m}^2}{16}\quad\textrm{by the fact}~w_1\geq\frac{1}{N_\mathfrak{m}}~
\textrm{and \eqref{yinequality}}.
\end{aligned}
\end{equation}
And if we let $\rho:=|q_{\omega_Q}|$, 
then we get by \eqref{yinequality} that
\begin{equation}\label{rhoinequality}
\rho=e^{-2\pi y_Q}\leq e^{-\pi\sqrt{3}}<\frac{1}{200}.
\end{equation}
Since $0\leq w_1\leq1/2$, we establish that
\begin{align}
\sum_{n=1}^\infty\left|\frac{q_{\omega_Q}^nq_{z_Q}}{(1-q_{\omega_Q}^nq_{z_Q})^2}\right|
&\leq\sum_{n=1}^\infty
\frac{\rho^{n+w_1}}{(1-\rho^{n+w_1})^2}
\leq\sum_{n=1}^\infty
\frac{\rho^{n}}{(1-\rho)^2}=\frac{\rho}{(1-\rho)^3},
\label{S3}
\\
\sum_{n=1}^\infty\left|\frac{q_{\omega_Q}^nq_{z_Q}^{-1}}{(1-q_{\omega_Q}^nq_{z_Q}^{-1})^2}\right|
&\leq\sum_{n=1}^\infty
\frac{\rho^{n-w_1}}{(1-\rho^{n-w_1})^2}
\leq\sum_{n=1}^\infty
\frac{\rho^{n-1/2}}{(1-\rho^{1/2})^2}=\frac{\rho^{1/2}}{(1-\rho)
(1-\rho^{1/2})^2},\label{S4}\\
2\sum_{n=1}^\infty\left|\frac{nq_{\omega_Q}^n}{1-q_{\omega_Q}^n}\right|
&\leq2\sum_{n=1}^\infty \frac{n\rho^n}{1-\rho^n}
\leq2\sum_{n=1}^\infty \frac{n\rho^n}{1-\rho}=\frac{2\rho}{(1-\rho)^3}.\label{S5}
\end{align}
We then derive that
\begin{align*}
\frac{1}{(2\pi)^2}|\wp(z_Q;\,[\omega_Q,\,1])|&<\frac{1}{12}+\frac{N_\mathfrak{m}^2}{16}
+\frac{\rho}{(1-\rho)^3}+\frac{\rho^{1/2}}{(1-\rho)(1-\rho^{1/2})^2}
+\frac{2\rho}{(1-\rho)^3}\\
&\hspace{7cm}\textrm{by \eqref{recallp} and \eqref{S20}--\eqref{S5}}\\
&=\frac{1}{12}+\frac{N_\mathfrak{m}^2}{16}
+\left(\frac{3\rho}{(1-\rho)^3}+\frac{\rho^{1/2}}{(1-\rho)(1-\rho^{1/2})^2}\right)\\
&<\frac{1}{12}+\frac{N_\mathfrak{m}^2}{16}+\frac{1}{10}\quad
\textrm{because}~\rho<\frac{1}{200}~\textrm{by \eqref{rhoinequality} and}\\
&\hspace{2cm}\frac{3x}{(1-x)^3}+\frac{x^{1/2}}{(1-x)(1-x^{1/2})^2}~\textrm{is increasing for}~0<x<1\\
&<\frac{N_\mathfrak{m}^2}{8}\quad\textrm{since}~N_\mathfrak{m}\geq2,
\end{align*}
from which it follows that
\begin{equation}\label{pinequality}
|\wp(z_Q;\,[\omega_Q,\,1])|<\frac{\pi^2}{2}N_\mathfrak{m}^2. 
\end{equation}
Using \eqref{g2expansion} in Lemma \ref{g2g3Deltaexpansion}, we deduce that
\begin{align*}
\frac{12}{(2\pi)^4}|g_2(\omega_Q)|&\leq1+240\sum_{n=1}^\infty\frac{n^3\rho^n}{1-\rho^n}\\
\\&<1+240\sum_{n=1}^\infty\frac{4^n\rho^n}{1-\rho}\quad\textrm{due to $n^3<4^n$ ($n\geq1$)}\\
&=1+240\frac{4\rho}{(1-\rho)(1-4\rho)}\quad\textrm{since $4\rho<1$ by \eqref{rhoinequality}}\\
&<6\quad\textrm{because $\rho<\frac{1}{200}$ by \eqref{rhoinequality}},
\end{align*}
and similarly we achieve that
\begin{align*}
\frac{6^3}{(2\pi)^6}|g_3(\omega_Q)|&\leq1+504\sum_{n=1}^\infty\frac{n^5\rho^n}{1-\rho^n}
\quad\textrm{by
\eqref{g3expansion}}\\
\\&<1+504\sum_{n=1}^\infty\frac{7^n\rho^n}{1-\rho}\quad\textrm{owing to the fact $n^5<7^n$ ($n\geq1$)}\\
&=1+504\frac{7\rho}{(1-\rho)(1-7\rho)}\quad\textrm{since $7\rho<1$ by \eqref{rhoinequality}}\\
&<20\quad\textrm{because $\rho<\frac{1}{200}$ by \eqref{rhoinequality}}.
\end{align*}
Furthermore, we obtain that
\begin{align*}
\frac{1}{(2\pi)^{12}}|\Delta(\omega_Q)|&\geq\rho\prod_{n=1}^\infty(1-\rho^n)^{24}\\
&\geq\rho\left(1-\sum_{n=1}^\infty\rho^n\right)^{24}\\
&=\rho\left(\frac{1-2\rho}{1-\rho}\right)^{24}\\
&>\rho\left(\frac{198}{199}\right)^{24}\quad\textrm{since $\rho<\frac{1}{200}$ by \eqref{rhoinequality}}\\
&>\frac{7}{8}e^{-2\pi y_Q}\\
&\geq\frac{7}{8}\,e^{-\pi\sqrt{|d_K|}}\quad
\textrm{because}~y_Q=\frac{\sqrt{|d_K|}}{2a}\leq\frac{\sqrt{|d_K|}}{2}.
\end{align*}
Thus it follows that
\begin{align}\label{g23deltainequality}
|g_2(\omega_Q)|<\frac{6(2\pi)^4}{12},\quad
|g_3(\omega_Q)|<\frac{20(2\pi)^6}{6^3},\quad
|\Delta(\omega_Q)|>\frac{7(2\pi)^{12}}{8}e^{-\pi\sqrt{|d_K|}}.
\end{align}
\begin{enumerate}
\item[Case 1.]
If $K\neq\mathbb{Q}(\sqrt{-1}),\,\mathbb{Q}(\sqrt{-3})$, then 
\begin{align*}
|f^{(k)}_\mathbf{w}(\omega_Q)|&=|f^{(1)}_\mathbf{w}(\omega_Q)|\\
&=|g_2(\omega_Q)|\,|g_3(\omega_Q)|\,\frac{1}{|\Delta(\omega_Q)|}\,|\wp(z_Q;\,[\omega_Q,\,1])|\\
&<\left(\frac{6(2\pi)^4}{12}\right)\left(\frac{20(2\pi)^6}{6^3}\right)
\left(\frac{8}{7(2\pi)^{12}}e^{\pi\sqrt{|d_K|}}\right)\left(\frac{\pi^2}{2}N_\mathfrak{m}^2\right)
\quad\textrm{by \eqref{pinequality} and \eqref{g23deltainequality}}\\
&=\frac{5}{756}\,N_\mathfrak{m}^2e^{\pi\sqrt{|d_K|}}\\
&<\mathcal{B}_\mathfrak{m}\quad\textrm{since}~N_\mathfrak{m}\geq2.
\end{align*}
\item[Case 2.] If $K=\mathbb{Q}(\sqrt{-1})$, then $\mathcal{Q}_\mathrm{red}(d_K)=\{x^2+y^2\}$ and so $\omega_Q=\mathrm{i}$. 
Note by Lemma \ref{g2g3Deltaexpansion} and the fact $0<q_\mathrm{i}=e^{-2\pi}<1$ that
\begin{equation}\label{g_2i}
\frac{12}{(2\pi)^4}|g_2(\mathrm{i})|>1.
\end{equation}
Then we have that
\begin{align*}
|f^{(k)}_\mathbf{w}(\omega_Q)|&=|f^{(2)}_\mathbf{w}(\mathrm{i})|\\
&=\frac{|\wp(w_1\mathrm{i}+w_2;\,[\mathrm{i},\,1])|^2}{|g_2(\mathrm{i})|}\quad\textrm{since $g_3(\mathrm{i})=0$ (see \cite[Theorem 3 in Chapter 3]{Lang})}\\
&<\left(\frac{\pi^2}{2}N_\mathfrak{m}^2\right)^2\frac{12}{(2\pi)^4}
\quad
\textrm{by \eqref{pinequality} and \eqref{g_2i}}\\
&<N_\mathfrak{m}^4\\
&<\mathcal{B}_\mathfrak{m}. 
\end{align*}
\item[Case 3.] If $K=\mathbb{Q}(\sqrt{-3})$, then $\mathcal{Q}_\mathrm{red}(d_K)=\{x^2+xy+y^2\}$ and so $\omega_Q=\zeta_3$.  
We see from \eqref{g23deltainequality} that
\begin{equation}\label{deltazeta_3}
|\Delta(\zeta_3)|>\frac{7(2\pi)^{12}}{8}e^{-\pi\sqrt{3}}. 
\end{equation}
And, we find that
\begin{align*}
|f^{(k)}_\mathbf{w}(\omega_Q)|&=|f^{(3)}_\mathbf{w}(\zeta_3)|\\
&=\frac{|\wp(w_1\zeta_3+w_2;\,[\zeta_3,\,1])|^3}{\sqrt{27|\Delta(\zeta_3)|}}\quad
\textrm{since $g_2(\zeta_3)=0$ (see \cite[Theorem 3 in Chapter 3]{Lang}) }\\
&\hspace{4.8cm}\textrm{and so $\Delta(\zeta_3)=
-27g_3(\zeta_3)^2$}\\
&<\left(\frac{\pi^2}{2}N_\mathfrak{m}^2\right)^3\sqrt{\frac{8e^{\pi\sqrt{3}}}{27\cdot7(2\pi)^{12}}}
\quad
\textrm{by \eqref{pinequality} and \eqref{deltazeta_3}}\\
&<N_\mathfrak{m}^6\\
&<\mathcal{B}_\mathfrak{m}. 
\end{align*}
\end{enumerate}
This completes the proof. 
\end{proof}

For a nonzero proper ideal  $\mathfrak{m}$ of $\mathcal{O}_K$, define an integer
\begin{equation}\label{Cm}
\mathcal{C}_\mathfrak{m}:=2+2^{3d_\mathfrak{m}}(\mathcal{A}_\mathfrak{m}\mathcal{B}_\mathfrak{m})^{2d_\mathfrak{m}}
\quad\textrm{where}~d_\mathfrak{m}:=[K_\mathfrak{m}:K].
\end{equation}

\begin{theorem}\label{intermediate}
Let $P\in E(\mathbb{C})$ be a primitive $\mathfrak{m}$-division point. 
If $K_\mathfrak{m}=K(h_E(P))$ (for example, when
$\mathfrak{f}_{K_\mathfrak{m}/K}=\mathfrak{m}$ by \textup{Theorem \ref{main}}), then
for every intermediate field $K\subseteq L\subseteq K_\mathfrak{m}$
\begin{equation}\label{relative}
L=K\bigg(
N_{K_\mathfrak{m}/L}(h_E(P)-\mathcal{C}_\mathfrak{m})
\bigg).
\end{equation}
\end{theorem}
\begin{proof}
Let $u=h_E(P)$. Since $K_\mathfrak{m}=K(u)$ by assumption, 
the Galois conjugates $u^g$ for $g\in\mathrm{Gal}(K_\mathfrak{m}/K)$
are all distinct. Furthermore, we note by Lemmas \ref{Aintegral} and \ref{hinequality} that
\begin{equation}\label{algebraic}
\textrm{$\mathcal{A}_\mathfrak{m}u$ is an algebraic integer}
\end{equation}
and
\begin{equation}\label{iotau}
|\iota(u)|\leq \mathcal{B}_\mathfrak{m}\quad\textrm{for every embedding}~\iota:K_\mathfrak{m}
\hookrightarrow\mathbb{C}. 
\end{equation} 
\begin{enumerate}
\item[Case 1.] If $L=K_\mathfrak{m}$, then 
\begin{equation*}
N_{K_\mathfrak{m}/L}(h_E(P)-\mathcal{C}_\mathfrak{m})=
h_E(P)-\mathcal{C}_\mathfrak{m}. 
\end{equation*}
Since $\mathcal{C}_\mathfrak{m}\in\mathbb{Z}$
and $K_\mathfrak{m}=K(h_E(P))$ by hypothesis, 
\eqref{relative} follows. 
\item[Case 2.] If $L=K$, then the relative norm belongs to $K$,
and so \eqref{relative} is obvious.
\item[Case 3.] Now, consider the case where $K\subsetneq L\subsetneq K_\mathfrak{m}$. If we put
\begin{equation*}
H=\mathrm{Gal}(K_\mathfrak{m}/L),\quad
d=|H|=[K_\mathfrak{m}:L],\quad \ell=[L:\mathbb{Q}],
\end{equation*}
then we observe that
\begin{equation}\label{degreeinequality}
d\geq2,\quad d\ell=[K_\mathfrak{m}:\mathbb{Q}]=2d_\mathfrak{m},\quad
\ell\leq d_\mathfrak{m}. 
\end{equation}
For $g\in\mathrm{Gal}(K_\mathfrak{m}/K)$, let
\begin{equation}\label{FgH}
F_{g,\,H}(x)=\prod_{\sigma\in H}(x-(u^{g})^\sigma)\in L[x].
\end{equation}
Furthermore, for $g\not\in H$, let
\begin{equation}\label{DgH}
D_{g,\,H}(x)=F_{g,\,H}(x)-F_{\mathrm{id},\,H}(x)\in L[x]. 
\end{equation}
Note that
if $g\not\in H$ and so $gH\cap H=\emptyset$, then the 
sets of roots of $F_{g,\,H}(x)$ and $F_{\mathrm{id},\,H}(x)$ are
disjoint because the associated Galois conjugates of $u$ are distinct. 
This implies that if $g\not\in H$, then the polynomial $D_{g,\,H}$ is 
nonzero. 
\par 
Now, let $g\not\in H$ and write
\begin{align*}
F_{g,\,H}(x)&=a_dx^d+a_{d-1}x^{d-1}+\cdots+a_0,\\
F_{\mathrm{id},\,H}(x)&=b_dx^d+b_{d-1}x^{d-1}+\cdots+b_0,\\
D_{g,\,H}(x)&=c_dx^d+c_{d-1}x^{d-1}+\cdots+c_0. 
\end{align*}
Then we achieve by the definitions \eqref{FgH}, \eqref{DgH} and \eqref{iotau} that 
for any embedding $\iota:L\hookrightarrow\mathbb{C}$ (which can be 
extended to an embedding of $K_\mathfrak{m}$ into $\mathbb{C}$) and 
$0\leq j\leq d$ 
\begin{align}\label{cab}
|\iota(c_j)|\leq|\iota(a_j)|+|\iota(b_j)|\leq \binom{d}{j}\mathcal{B}_\mathfrak{m}^{d-j}+
\binom{d}{j}\mathcal{B}_\mathfrak{m}^{d-j}\leq U,\quad
\textrm{where}~U=2^{d+1}\mathcal{B}_\mathfrak{m}^d.
\end{align}
\par
And, we know by \eqref{algebraic} that the coefficients of
$\mathcal{A}_\mathfrak{m}^dF_{g,\,H}(x)$ and $\mathcal{A}_\mathfrak{m}^dF_{\mathrm{id},\,H}(x)$
are all algebraic integers. 
Let $\delta$ be the nonzero leading coefficient of $D_{g,\,H}(x)$. 
We then derive that
\begin{align*}
1&\leq|N_{L/\mathbb{Q}}(\mathcal{A}_\mathfrak{m}^d\delta)|\quad
\textrm{since $\mathcal{A}_\mathfrak{m}^d\delta$ is a nonzero algebraic integer}\\
&=\mathcal{A}_\mathfrak{m}^{d\ell}|N_{L/\mathbb{Q}}(\delta)|\\
&=\mathcal{A}_\mathfrak{m}^{d\ell}|\delta|\cdot\prod_\iota|\iota(\delta)|
\quad\textrm{where $\iota$ runs over all nonidentity embeddings $L\hookrightarrow\mathbb{C}$}\\
&\leq \mathcal{A}_\mathfrak{m}^{d\ell}|\delta|U^{\ell-1}\quad\textrm{by \eqref{cab}},
\end{align*}
and hence
\begin{equation}\label{deltainequality}
|\delta|\geq \mathcal{A}_\mathfrak{m}^{-d\ell}U^{-(\ell-1)}.
\end{equation}
If $D_{g,\,H}(x)$ is constant, then it has no zeros. Otherwise, 
every zero $z$ of $D_{g,\,H}(x)$ satisfies
\begin{align*}
|z|&\leq 1+\frac{U}{|\delta|}\quad\textrm{by Cauchy's bound (see
\cite[p. 123]{Marden})
and \eqref{cab}}\\
&\leq 1+\mathcal{A}_\mathfrak{m}^{d\ell}U^\ell\quad\textrm{by \eqref{deltainequality}}\\
&=1+2^{2d_\mathfrak{m}+\ell}(\mathcal{A}_\mathfrak{m}\mathcal{B}_\mathfrak{m})^{2d_\mathfrak{m}}
\quad\textrm{because $d\ell=2d_\mathfrak{m}$ by \eqref{degreeinequality} and $U=2^{d+1}\mathcal{B}_\mathfrak{m}^d$}\\
&\leq1+2^{3d_\mathfrak{m}}(\mathcal{A}_\mathfrak{m}\mathcal{B}_\mathfrak{m})^{2d_\mathfrak{m}}\quad
\textrm{since}~\ell\leq d_\mathfrak{m}~\textrm{by \eqref{degreeinequality}}\\
&<\mathcal{C}_\mathfrak{m}\quad\textrm{by the definition \eqref{Cm}}. 
\end{align*}
This observation yields that 
\begin{equation}\label{Dnot0}
D_{g,\,H}(\mathcal{C}_\mathfrak{m})\neq0\quad\textrm{for every}~g\in\mathrm{Gal}(K_\mathfrak{m}/K)
\setminus H. 
\end{equation}
\par
Let
\begin{equation*}
v=N_{K_\mathfrak{m}/L}(u-\mathcal{C}_\mathfrak{m}).
\end{equation*}
We find by the definition \eqref{FgH} that for $g\in\mathrm{Gal}(K_\mathfrak{m}/K)$ 
\begin{equation*}
v^g=(-1)^dF_{g,\,H}(\mathcal{C}_\mathfrak{m}).
\end{equation*}
Thus every element of $H$ fixes $v$, whereas no element of $\mathrm{Gal}(K_\mathfrak{m}/K)
\setminus H$ fixes $v$ by the definition \eqref{DgH} and \eqref{Dnot0}. 
Therefore we conclude by Galois theory that
\begin{equation*}
L=K(v)=K\bigg(
N_{K_\mathfrak{m}/L}(h_E(P)-\mathcal{C}_\mathfrak{m})
\bigg).
\end{equation*}
\end{enumerate}
\end{proof}

\section {Hilbert's twelfth problem for imaginary quadratic fields}\label{sect9}

We combine the relative norm construction with the description 
of ray class groups by extended form class groups. 
This gives explicit generators for all finite abelian extensions of an 
imaginary quadratic field in terms of special values of Fricke functions.

\begin{lemma}\label{4lemma}
If $\mathfrak{m}=N\mathcal{O}_K$ for a positive integer $N$ divisible by $4$, then we have
\begin{equation*}
K_\mathfrak{m}=K\bigg(h_E(\phi_E(\tfrac{1}{N}+\mathcal{O}_K))\bigg).
\end{equation*}
\end{lemma}
\begin{proof}
Note that the point $\phi_E(\frac{1}{N}+\mathcal{O}_K)\in E(\mathbb{C})$ is a primitive 
$\mathfrak{m}$-division point. 
If $N=4$, the assertion follows from \cite[Corollary 7.1]{J-K-S}. 
If $N>4$, then $\mathfrak{m}$ satisfies (S1), so 
the assertion follows from \cite[p. 190]{Sugawara1}. 
\end{proof}

For a positive integer $N$, let $\mathcal{Q}(d_K,\,N)$
be the subset of $\mathcal{Q}(d_K)$ defined by
\begin{equation*}
\mathcal{Q}(d_K,\,N)=\{ax^2+bxy+cy^2\in\mathcal{Q}(d_K)\mid
\gcd(a,\,N)=1\}. 
\end{equation*} 
The congruence subgroup
\begin{equation*}
\Gamma_1(N)=\left\{\gamma\in\mathrm{SL}_2(\mathbb{Z})\mid
\gamma\equiv\begin{bmatrix}1&*\\0&1\end{bmatrix}\Mod{NM_2(\mathbb{Z})}\right\}
\end{equation*}
acts from the right on the set $\mathcal{Q}(d_K,\,N)$ by 
\begin{equation*}
Q^\gamma =Q\left(\gamma\begin{bmatrix}x\\y\end{bmatrix}\right)\quad
(Q\in\mathcal{Q}(d_K,\,N),~\gamma\in\Gamma_1(N)),
\end{equation*}
which yields an equivalence relation $\sim_{\Gamma_1(N)}$ on $\mathcal{Q}(d_K,\,N)$. 
For $Q\in\mathcal{Q}(d_K,\,N)$, let $[Q]_{\Gamma_1(N)}$
be its equivalence class under $\sim_{\Gamma_1(N)}$.
Furthermore, let $C_{\Gamma_1(N)}(d_K,\,N)$ denote the set of these equivalence classes, namely,
\begin{equation*}
C_{\Gamma_1(N)}(d_K,\,N)=\mathcal{Q}(d_K,\,N)/\sim_{\Gamma_1(N)}
=\{[Q]_{\Gamma_1(N)}\mid Q\in\mathcal{Q}(d_K,\,N)\}.  
\end{equation*}
For $Q\in\mathcal{Q}(d_K,\,N)$, 
write 
\begin{equation*}
Q=a_Qx^2+b_Qxy+c_Qy^2,
\end{equation*}
and choose $a_Q'\in\mathbb{Z}$ with $a_Qa_Q'\equiv1\Mod{N}$. 
Recall that $\omega_Q$ is the zero of $Q(x,\,1)$ in $\mathbb{H}$. 
Put
\begin{equation*}
b_K=\left\{\begin{array}{ll}0 & \textrm{if}~d_K\equiv0\Mod{4},\\
1 & \textrm{if}~d_K\equiv1\Mod{4}\end{array}\right.
\end{equation*} 
so that $\tau_K=(-b_K+\sqrt{d_K})/2$. 

\begin{proposition}\label{extended}
Let $\mathfrak{m}=N\mathcal{O}_K$ for a positive integer $N$. 
There is a binary operation on $C_{\Gamma_1(N)}(d_K,\,N)$
which makes $C_{\Gamma_1(N)}(d_K,\,N)$ a group isomorphic to 
the ray class group $C(\mathfrak{m})$ via the mapping
\begin{equation*}
C_{\Gamma_1(N)}(d_K,\,N)\rightarrow
C(\mathfrak{m}),\quad
[Q]_{\Gamma_1(N)}\mapsto[[\omega_Q,\,1]]_\mathfrak{m}. 
\end{equation*}
Moreover, the mapping
\begin{eqnarray*}
\sigma_N:C_{\Gamma_1(N)}(d_K,\,N)&\rightarrow&
\mathrm{Gal}(K_\mathfrak{m}/K),\\
\mathrm{[}Q]_{\Gamma_1(N)}&\mapsto&
\bigg(
f(\tau_K)\mapsto f^{\left[\begin{smallmatrix}
1 & -a_Q'\left(\frac{b_Q+b_K}{2}\right)\\0&a_Q'
\end{smallmatrix}\right]}(-\overline{\omega}_Q)\mid f\in\mathcal{F}_N~\textrm{is finite at}~\tau_K\bigg)
\end{eqnarray*}
is a well-defined isomorphism.
\end{proposition}
\begin{proof}
See \cite[Theorems 9.4 and 12.3]{J-K-S-Y}. 
\end{proof}

\begin{remark}
We call the group $C_{\Gamma_1(N)}(d_K,\,N)$
the \textit{extended form class group of discriminant $d_K$ for $\Gamma_1(N)$}.
When $N=1$, this group is none other than 
Gauss's classical form class group (see \cite[Theorems 3.9 and 7.7]{Cox}). 
\end{remark}

\begin{theorem}\label{Hilbert}
Let $K$ be an imaginary quadratic field of discriminant $d_K$, 
and let $L$ be a finite extension of $K$. 
Then $L/K$ is abelian if and only if 
there exist a positive integer $N$ divisible by $4$ and a subgroup $H$ of 
the extended form class group $C_{\Gamma_1(N)}(d_K,\,N)$
such that
\begin{equation*}
L=K\left(
\prod_{[Q]_{\Gamma_1(N)}\in H}\left(f_{\left[\begin{smallmatrix}
0 &a_Q'/N
\end{smallmatrix}\right]}^{(|\mathcal{O}_K^\times|/2)}(-\overline{\omega}_Q)
-\mathcal{C}_{N\mathcal{O}_K}\right)\right).
\end{equation*}
\end{theorem}
\begin{proof}
For a positive integer $N$, let $\sigma_N:C_{\Gamma_1(N)}(d_K,\,N)\stackrel{\sim}{\rightarrow}
\mathrm{Gal}(K_{N\mathcal{O}_K}/K)$ be the isomorphism stated in Proposition \ref{extended}. 
Then we establish that
\begin{eqnarray*}
&&\textrm{$L/K$ is abelian}\\
&\Longleftrightarrow&L\subset K_\mathfrak{m}~\textrm{for some nonzero ideal
$\mathfrak{m}$ of $\mathcal{O}_K$}\\
&\Longleftrightarrow&L\subseteq K_{N\mathcal{O}_K}~\textrm{for some 
positive integer $N$ divisible by $4$}\\
&\Longleftrightarrow&\textrm{$L$ is the 
subfield of $K_{N\mathcal{O}_K}$ which 
is the fixed field of a subgroup of $\mathrm{Gal}(K_{N\mathcal{O}_K}/K)$}\\
&&\textrm{for some positive integer $N$ divisible by $4$}\\
&\Longleftrightarrow&\textrm{$L$ is the subfield
of $K_{N\mathcal{O}_K}$ which is the fixed field
of $\sigma_N(H)$}\\
&&\textrm{with $H$ a subgroup of $C_{\Gamma_1(N)}(d_K,\,N)$ 
for some positive integer $N$ divisible by $4$},\\
&&\textrm{by Proposition \ref{extended}}\\
&\Longleftrightarrow&L=K\left(
\prod_{[Q]_{\Gamma_1(N)}\in H}\bigg(f_{\left[\begin{smallmatrix}0&1/N\end{smallmatrix}\right]}
^{(|\mathcal{O}_K^\times|/2)}(\tau_K)-\mathcal{C}_{N\mathcal{O}_K}\bigg)^{\sigma_N([Q]_{\Gamma_1(N)})}
\right)\\
&&\textrm{with $H$ a subgroup of $C_{\Gamma_1(N)}(d_K,\,N)$ for some positive integer $N$ divisible by $4$,}
\\&&\textrm{by Theorem \ref{intermediate}, Lemma \ref{4lemma},
the definitions \eqref{Weberh} and \eqref{Frickedef}}\\
&\Longleftrightarrow&L=K\left(
\prod_{[Q]_{\Gamma_1(N)}\in H}\left(f_{\left[\begin{smallmatrix}
0 &a_Q'/N
\end{smallmatrix}\right]}^{(|\mathcal{O}_K^\times|/2)}(-\overline{\omega}_Q)
-\mathcal{C}_{N\mathcal{O}_K}\right)\right)\\
&&\textrm{with $H$ a subgroup of $C_{\Gamma_1(N)}(d_K,\,N)$ for some positive integer $N$ divisible by $4$},
\\&&\textrm{by Propositions \ref{modularGalois} (ii), \ref{extended} and \eqref{f23}.}
\end{eqnarray*}
\end{proof}

\section*{Statements \& Declarations}

\subsection*{Funding}
The first author was supported
by the National Research Foundation of Korea (NRF) grant funded by the Korea government (MSIT)
(No. RS-2026-25606107).
The second author was supported
by the National Research Foundation of Korea (NRF) grant funded by the Korea government (MSIT)
(No. RS-2026-25588491).
The third author was supported
by the National Research Foundation of Korea (NRF) grant funded by the Korea government (MSIT)
(No. RS-2026-25588491) and
Hankuk University of Foreign Studies Research Fund of 2026.

\subsection*{Competing interests}
The authors have not disclosed any competing interests.

\subsection*{Data availability}
Data sharing is not applicable to this article as no datasets
were generated or analysed during the current study.

\bibliographystyle{amsplain}

\address{
Department of Mathematics\\
Dankook University\\
Cheonan-si, Chungnam 31116\\
Republic of Korea} {hoyunjung@dankook.ac.kr}
\address{
Department of Mathematical Sciences \\
KAIST \\
Daejeon 34141\\
Republic of Korea} {jkgoo@kaist.ac.kr}
\address{
Department of Mathematics\\
Hankuk University of Foreign Studies\\
Yongin-si, Gyeonggi-do 17035\\
Republic of Korea} {dhshin@hufs.ac.kr}

\end{document}